\documentclass[10pt]{article}

\usepackage{amsfonts} \usepackage{latexsym}

 \RequirePackage{amsbsy} % for Boldsymbol
 \RequirePackage{amsopn} % for DeclareMathOperator
 \RequirePackage{amsmath}% for multiple lines in equations
  \RequirePackage{amssymb}% for subsetneq

 \newtheorem{thee}{Theorem}[section] \newtheorem{coor}[thee]{Corollary}
 \newtheorem{leem}[thee]{Lemma} \newtheorem{prro}[thee]{Proposition}
  \newtheorem{exxe}[thee]{Example}
 \newtheorem{reem}[thee]{Remark}

 \newcommand{\balf} {\renewcommand{\theenumi}{(\alph{enumi})}
 \renewcommand{\labelenumi}{\theenumi}
                      \begin{enumerate}}
\newcommand{\ealf}   {\end{enumerate}
                      \renewcommand{\theenumi}{\arabic{enumi}}
                      \renewcommand{\labelenumi}{\theenumi.}}
\newcommand{\bara}   {\renewcommand{\theenumi}{(\arabic{enumi})}
                      \renewcommand{\labelenumi}{\theenumi}
                      \begin{enumerate} }
\newcommand{\eara}   {\end{enumerate}
                      \renewcommand{\theenumi}{\arabic{enumi}}
                      \renewcommand{\labelenumi}{\theenumi.}}
 
 \newcommand{\brom}   {\renewcommand{\theenumi}{(\roman{enumi})}
                      \renewcommand{\labelenumi}{\theenumi}
                      \begin{enumerate} }
\newcommand{\erom}   {\end{enumerate}
                      \renewcommand{\theenumi}{\arabic{enumi}}
                      \renewcommand{\labelenumi}{\theenumi.}}

 \DeclareMathOperator{\Na}{Na}

\begin{document}

 \title {\bf Semistar linkedness and flatness, \\
 Pr\"ufer semistar multiplication domains} \medskip
 
 \author{\bf Said El Baghdadi$^{1}$\, \rm and \, \bf Marco Fontana$^{2}$
 \footnote{During the preparation of this work, the second named 
 author was 
 supported in part by a research grant MIUR 2001/2002 (Cofin
 2000-MM01192794).}
    \\
 \rm \\
 \small $^{1}$ Department of Mathematics, Facult\'e des Sciences et
 Techniques, \\
 \small P.O. Box 523, Beni Mellal, Morocco --  \texttt{\small
 baghdadi@fstbm.ac.ma}\\
 \small  $^{2}$ Dipartimento di Matematica, Universit\`a degli Studi  ``Roma
 Tre'', \\
 \small  Largo S. L. Murialdo, 1 - 00146 Roma, Italy --
\small \texttt{\small fontana@mat.uniroma3.it}\\}

 \date{ }
 
 \maketitle
 
 \bigskip 
 
 \begin{abstract}
     In 1994, Matsuda and Okabe introduced the notion of 
semistar operation, extending the ``classical" concept of star operation.  In
this paper, we introduce and study the notions of semistar linkedness
and semistar flatness which are natural generalizations, to the semistar
setting, of their corresponding ``classical" concepts.  As an
application, among other results, we obtain a semistar version of Davis'
and Richman's overring-theoretical theorems of characterization of
Pr\"ufer domains for Pr\"ufer semistar multiplication domains.
     \end{abstract}
     
     \bigskip \bigskip
 
%%%%%%%%%%%%%%%%%%%%%%%%%%%%%%%%%%%%%%%%%%%%%%%%%%%%%%%%%%%%%%%%%%%%%
\small
    
      \section{{Introduction}}\label{sec:1}

 \hskip 0.5cm   Star operations have a central place in multiplicative ideal theory, this 
concept arises from the classical theory of ideal systems, based on the 
work by W. Krull, E. Noether, H. Pr\"ufer, and P. Lorenzen (cf. 
\cite{GILMER:1972}, \cite{Jaffard}, \cite{HK}).  Recently, new interest
on these theories has been originated by the work by R. Matsuda and A. Okabe
\cite{Okabe/Matsuda:1994}, where the notion of
semistar operation was introduced, as a generalization of the notion of star operation. 
This concept has been proven, regarding its flexibility, extremely
useful in studying the structure of different classes of integral domains
(cf.  for instance \cite{Matsuda/Sugatani:1995}, \cite{FH}, \cite{FL1},
\cite{FL2}, \cite{FL3}, and \cite{Koch:JA}).

\hskip 0.5cm Recall that a domain $D$, on which a semistar operation
$\star$ is defined, is called a \it Pr\"ufer semistar multiplication
domain \rm (or \it P$\star$MD\rm ), if each nonzero finitely generated ideal $F$
of $D$ is $\star_f$-invertible (i.e., $(FF^{-1})^{\star_f}=D^{\star}$),
where $\star_f$ is the semistar operation of finite type associated to
$\star$ (cf.  Section 2 for details).  These domains generalize Pr\"ufer
$v$--multiplication domains \cite[page 427]{GILMER:1972} (and, in
particular, Pr\"ufer and Krull domains) to the semistar multiplication setting.

\vskip -3pt Among the various overring-theoretical characterization of Pr\"ufer
domains, the following two have relevant consequences:

$\bullet$\; Davis' characterization \cite[Theorem 1]{Davis}: \sl a domain
$D$ is a Pr\"ufer domain if and only if each overring of $D$ is
integrally closed; \rm

$\bullet$\; Richman's characterization \cite[Theorem 4 ]{R}: \sl a domain $D$ is a
Pr\"ufer domain if and only if each overring of $D$ is $D$--flat.  \rm

The previous theorems have been extended to the case of Pr\"ufer
$v$--multiplication domains (for short, P$v$MDs) in \cite{DHLZ1} and \cite{KP},
respectively, by means of the $v$ (or the $t$)--operation.

\hskip 0.5cm The purpose of the present work is to deepen the study of a
general multiplicative theory for the semistar context, with special
emphasis to the linkedness and the flatness, and to pursue the study of
Pr\"ufer semistar multiplication domains (cf. 
\cite{Houston/Malik/Mott:1984}, \cite{FJS}).

\vskip -3pt
In Section 2 we recall the main definitions and we collect some background
results on semistar operations.  In Section 3, we define and study the
notion of semistar linked overring, which generalizes the notion of
$t$--linked overring defined in \cite{DHLZ1}.  Several characterizations of this 
concept have been obtained.  Section 4 is devoted to semistar flat
overrings, a concept which extends the classical notion of flat overring
and gives a very ``flexible" general tool, preserving for the ``semistar
prime ideals" involved,  a similar behaviour as in the classical context.  As an
application, in Section 5, we achieve the proofs for analogues of Davis' and
Richman's theorems in the general case of Pr\"ufer semistar multiplication
domains.  \bigskip \bigskip

 %\begin{center}
       \section{Background and preliminary results on semistar operations}
   %\end{center}

      \vskip -9pt \hskip 0.5cm Let $D$ be an integral domain with quotient
      field $K$.  Let $\boldsymbol{\overline{F}}(D)$ denote the set of all
      nonzero $D$-submodules of $K$ and let $\boldsymbol{F}(D)$ be the set of
      all nonzero fractional ideals of $D$, i.e. all $E \in
      \boldsymbol{\overline{F}}(D) $ such that there exists a nonzero $d \in
      D$ with $dE \subseteq D$.  Let $\boldsymbol{f}(D)$ be the set of all
      nonzero finitely generated $D$-submodules of $K$.  Then, obviously $\,
      \boldsymbol{f}(D) \subseteq \boldsymbol{F}(D) \subseteq
      \boldsymbol{\overline{F}}(D) \, .$

\hskip 0.5cm We recall that a mapping
$ \star : \boldsymbol{\overline{F}}(D) \rightarrow
\boldsymbol{\overline{F}}(D) \,, \;\;\;E \mapsto E^{\star} $
 is called a \it semistar operation on $D$ \rm if, for $x
\in K, x \not = 0$, and $E,F \in \boldsymbol{\overline{F}}(D) $, the
following properties hold:

\hspace*{30pt}$\boldsymbol{(\star_1)}\;$ $(xE)^{\star} = xE^{\star}\,; $

\hspace*{30pt}$\boldsymbol{(\star_2)}\;$  $E \subseteq F
\Rightarrow  E^{\star} \subseteq F^{\star}\,; $

\hspace*{30pt}$\boldsymbol{(\star_3)}\;$  $E \subseteq E^{\star}$ \,  and \, 
$E^{\star} =
(E^{\star})^{\star} =: E^{\star \star}$\,,

\noindent cf.  for instance \cite{OM1}, \cite{Okabe/Matsuda:1994},
\cite{Matsuda/Sugatani:1995}, \cite{Matsuda/Sato:1996}, \cite{FH} and
\cite{FL1}. 

\vskip -3 pt When $D^\star=D$, the semistar operation $\star$, restricted
to $\boldsymbol{{{F}}}(D)$, is ``the classical'' star operation (cf.  \cite
[Sections 32 and 34]{GILMER:1972}).  In this case, we will write that $\star$ is
\it a (semi)star operation on $D$\rm .  \par

%EXAMPLE 2.1
\begin{exxe} \label{ex:2.1} \bf (1) \rm The constant map $E \mapsto E^e:=K$,
$E\in \boldsymbol{{\overline{F}}}(D)$, defines a trivial semistar
operation $e$ (or, $e_{D}$) on $D$, called \it  the $e$-operation\rm .

\bf (2) \rm The map $E \mapsto E^d:=E$,
$E\in\boldsymbol{{\overline{F}}}(D)$, defines a (semi)star operation $d$
(or, $d_{D}$) on $D$, called \it  the $d$--operation \rm or \it the identity
semistar operation\rm .

\bf (3) \rm For each $E\in\boldsymbol{{\overline{F}}}(D)$, set $E^{-1}:=(D:_KE):=\{x\in K,\,
xE\subseteq D\}$.  The map $E \mapsto E^v:=(E^{-1})^{-1}$ defines a (semi)star
operation on $D$, called \it  the $v$--operation \rm on $D$ (or \it  the
$v_{D}$--operation\rm ).  This operation, when restricted to
$\boldsymbol{{ {F}}}(D)$, is the classical $v$--operation on $D$.

\bf (4) \rm Let $\{T_\lambda \mid \lambda \in \Lambda \}$ be a
family of overrings of $D$, and let $\ast_{\lambda}$ be a semistar
operation on $T_{\lambda}$, for each $\lambda \in \Lambda$. Then $E \mapsto
 E^{\ast_{\Lambda}}
:= \cap \{ (ET_\lambda)^{\ast_{\lambda}} \mid \lambda \in \Lambda \}$, is a
semistar operation on $D$.  Moreover, $(E^{\ast_{\Lambda}} T_\lambda)^{\ast_{\lambda}}
=(ET_\lambda)^{\ast_{\lambda}}$, for each $\lambda \in \Lambda$.  This semistar
operation is called \it the semistar operation induced by the family \rm 
$\{(T_\lambda, \ast_{\lambda}) \mid \lambda \in \Lambda \}$ (for the
star case, cf.  \cite[Theorem 2]{Anderson:1988} and, for the semistar
case, cf.  \cite[Example 1.3 (d)]{FH},  \cite[Example 2.1 (g)]{FJS}).  Note
that, in general, $D$ is a proper subset of $D^{\ast_{\Lambda}}= \cap \{
(T_\lambda)^{\ast_{\lambda}} \mid \lambda \in \Lambda \}$.  In
particular, if $T$ is an overring of $D$, we denote by $\star_{\{T\}}$
the semistar operation induced by $\{(T, d_{T})\}$.  For example, we
have that $e_{D}=\star_{\{K\}}$ and $d_{D}=\star_{\{D\}}$.

\bf (5) \rm Spectral semistar operations constitute perhaps the most
important class of semistar operations of the type introduced in (4). 
Given a set $\Delta$ of prime ideals of an integral domain $D$, \it the spectral semistar operation $\star_{\Delta}$  on $D$
  associated to $\Delta$ \rm is the semistar operation on $D$  induced by the family $\{(D_P, d_{D_{P}})
 \mid \, P\in \Delta\}$ (cf.  the previous Example (4)); when $\Delta=
 \emptyset$, then we set $\star _{\emptyset}:=e_{D}$. A \it spectral semistar
 operation on $D$ \rm is a semistar operation $\star$ on $D$ such that
 there exists a set of prime ideals $\Delta$ of $D$ with $\star = \star
 _{\Delta}$.  A spectral semistar operation $\star$ is \it a stable
 semistar operation\rm , i.e., $(E\cap F)^\star=(E^\star\cap F^\star)$,
 for all $E, F\in{\overline{\boldsymbol{F}}}(D)$ (or, equivalently, $(E:
 F)^\star=(E^\star: F^\star)$, for each $E\in
 {\overline{\boldsymbol{F}}}(D)$ and $F\in \boldsymbol{{f}}(D))$.  For
 more details, see \cite[Section 4]{FH}.

\bf (6) \rm Let $D$ be an integral domain and $T$ an
overring of $D$.  Let $\star$ be a semistar operation on $D$, the map
$\dot{\star}^{T}$: ${\overline{\boldsymbol{F}}}(T) \rightarrow {\overline{\boldsymbol{F}}}(T)$,
$E^{\dot{\star}^{T}}:=E^\star$, for $E\in {\overline{\boldsymbol{F}}}(T)\ (\subseteq
{\overline{\boldsymbol{F}}}(D)) $, is a semistar operation on $T$.  When
$T:=D^\star$, then we set simply $\dot{\star}$, instead of 
$\dot{\star}^{D^\star}$, and we note that $\dot{\star}$ is a (semi)star
operation on $D^\star$.

\vskip -3pt Conversely, let $\ast$ be a semistar operation on an overring $T$ of
$D$ and define $\hbox{\d{$\ast$}}_{D}$: ${\overline{\boldsymbol{F}}}(D)
\rightarrow {\overline{\boldsymbol{F}}}(D)$, by setting
$E^{\hbox{\d{$\ast$}}_{D}}:=(ET)^\ast$, for each $E\in
{\overline{\boldsymbol{F}}}(D)$.  For each semistar operation $\ast$ on $T$, if
we set $\star:= {\hbox{\d{$\ast$}}}_{D}$, then we have that $
{\dot{\star}}^{T}=\ast$ \cite[Corollary 2.10]{FL1}.

\bf (7) \rm  Given a semistar operation $\star$ on $D$, we can
define a new semistar operation on $D$, by setting $E \mapsto E^{\star_f}
:=\cup\{F^\star \mid \, F\in {\boldsymbol{f}}(D),\, F\subseteq E\}$, for each 
$E\in {\overline{{\boldsymbol{F}}}}(D)$.  The semistar operation $\star_f$ is
called \it the semistar operation of finite type \rm associated to $\star$.   Note 
that if $E\in {\boldsymbol{f}}(D)$, then $E^\star=E^{\star_f}$. 
A semistar operation $\star$ is called \it a semistar operation of finite
type \rm if $\star=\star_f$.  Note that $(\star_{f})_{f} = \star_{f}$.

\vskip -3pt An important example of semistar operation of finite type is the
(semi)star operation of finite type associated to the $v$-(semi)star
operation, i.e. $t:=v_f$, called \it the $t$--(semi)star operation on
$D$ \rm (or \it the $t_{D}$--operation\rm ).  Note that the $e$--operation and the
identity operation $d$ on $D$ are of finite type.  A spectral semistar
operation $\star $ on $D$ is of finite type if and only if $\star=
\star_{\Delta}$, for some quasi--compact set of prime ideals of $D$
\cite[Corollary 4.6 (2)]{FH}.

\end{exxe}

 \hskip .5cm If $\star_1$ and $\star_2$ are two semistar operations on
an integral domain  $D$, we say that $\, \star_1\le\star_2\,$ if , 
for each $E\in {\overline{{\boldsymbol{F}}}}(D)$, $E^{\star_1}\subseteq
E^{\star_2}$; in this case
$(E^{\star_1})^{\star_2}=E^{\star_2}$.

\vskip -3pt Note that, for each semistar operation $\star$, we have
that $\star_f \le \star$.  Moreover, for each (semi)star operation $\star$ on $D$,
we have always that $\star \le v$ and, hence, $\star_f \le t$
(easy consequence of \cite[Theorem 34.1 (4)]{GILMER:1972}).

	\hskip 0.5cm \rm Let $\,I \subseteq D\,$ be a nonzero ideal of $\,D\,$
	and let $\,\star\,$ be a semistar operation on $\,D\,$.  We say that
	$\,I\,$ is \it a quasi--$\star$--ideal \rm (respectively,
	\it $\star$--ideal\rm ) \it of $\,D\,$ \rm if $\,I^\star \cap D = I\,$
	(respectively, $\,I^\star = I\,$).  \ Similarly, we call \it a
	quasi--$\star$--prime \rm (respectively,\, \it a $\,\star$--prime \rm)
	of $D$ a quasi--$\star$--ideal (respectively,\, $\,\star$--ideal) of
	$\,D\,$ which is also a prime ideal.  We call \it a
	quasi--$\star$--maximal \rm (respectively,\, \it a $\,\star$--maximal
	\rm) of $\,D\,$ a maximal element in the set of all proper
	quasi--$\star$--ideals (respectively,\, $\star$--ideals) of
	$\,D\,.$
	
	\hskip 0.5cm  Note that if $\,I \subseteq D\,$ is a $\,\star$--ideal, it is also a
quasi--$\star$--ideal and, when $\,D = D^\star\,$,
the notions of quasi--$\star$--ideal and
 $\,\star$--ideal coincide.
\smallskip

 \hskip 0.5cm When $\,D \subsetneq D^\star \subsetneq K\,$ we can ``restrict'' the
 semistar operation $\,\star\,$ on $\,D\,$ to the (semi)star operation
 $\,{\dot{\star}}\,$ on $\,{D^\star}\,$ (Example \ref{ex:2.1} (6)) and we have a
 strict relation between the quasi--$\star$--ideals of $\,D\,$ and the
 $\,{\dot{\star}}$--ideals of $\,{D^\star}\,$, \, as shown in the following
 result:

%LEMMA 2.2
\begin{leem} \label{lm:2.2} \rm \cite [Lemma 2.2]{FL3}.  \sl Let
$D$ be an integral domain and $\star$ a semistar operation on $D$ and
let $\,{\dot{\star}}\,$ be the (semi)star operation on $\,{D^\star}\,$
associated to $\star$.  Then: 

\balf \rm \bf \item \sl $ I \, \hbox{ is a
quasi--$\star$--ideal of $D$} \Leftrightarrow I = L \cap D,  \,
\hbox{where} \, L \subseteq D^{\star} \, \hbox{ is a
$\dot{\star}$--ideal of $D^\star$}.$

\rm \bf \item \sl If $\, L \subseteq D^{\star} \,$ is a $\dot{\star}$--prime
ideal of $\,D^\star,$ then $\,L\cap D\,$ is a quasi--$\star$--prime ideal of
$\,D.$ \hfill $\Box$ \ealf
\end{leem}

 \hskip 0.5cm   Note that, in general, the restriction to $\,D\,$ of a
$\dot{\star}$--maximal ideal of $D^\star$  is a
quasi--$\star$--prime ideal of $\,D\,,$  but not necessarily a
quasi--$\star$--maximal ideal of $\,D\,,$  and if $\, L\,$ is an 
ideal of
$\,D^\star\,$ and $\,L\cap D\,$ is a quasi--$\star$--prime ideal of
$\,D\,,$ then $\, L\,$ is not necessarily a $\dot{\star}$--prime ideal
of $\,D^\star$, \cite [Remark 3.6]{FL3}.

%LEMMA 2.3
\begin{leem} \label{lm:2.3} \sl Let $\,\star\,$ be a
semistar operation of an integral domain $\,D\,$.  Assume that $\,\star\,$ is
not trivial and that $\,\star = \star_f\,$.  Then:
\balf

\rm \bf \item \sl Each proper
quasi--$\star$--ideal is contained in a quasi-$\star$-maximal.

\rm \bf \item \sl Each
quasi--$\star$--maximal is a quasi--$\star$-prime.

\rm \bf \item \sl If $\, Q\, \hbox{ is a quasi--$\star$--maximal ideal of $D$}\,$
then $\, Q = M \cap D\,, \, $ for some $\,\dot{\star}$--maximal ideal $\,M\,$ of
$D^\star.$

\rm \bf \item \sl Each minimal prime over a
quasi--$\star$--ideal is a quasi--$\star$--prime.

\rm \bf \item \sl  Set
$$
\Pi^\star : = \{P \in \hbox{\rm Spec}(D)\, \mid\,\, P \not = 0 \, 
\hbox{
and } \, P^\star \cap D \not = D\}\,,
$$
\noindent then each quasi--$\star$-prime of $D$ belongs to $\Pi^\star$ and,
moreover, the set of maximal elements of $\,\Pi^\star$ is nonempty and
coincides with the set of all the quasi--$\star$-maximals of $D$. 
\ealf
\end{leem}
\bf Proof.  \rm We give a proof of \bf (d)\rm , for the other statements
see \cite [Lemma 2.3]{FL3}.

Let $I$ be a  quasi--$\star$--ideal of $D$ and let $P$ a minimal prime
ideal of $D$ over $I$, hence rad$(ID_{P}) =PD_{P}$.  Then, for each
finitely generated ideal $J$ of $D$, with $J \subseteq P$, there exists
an integer $n \geq 1$ such that $J^nD_{P}\subseteq ID_{P}$, i.e.
$sJ^n\subseteq I$, for some $s \in D\setminus P$.  Therefore:
$$
\begin{array}{rl}
& s(J^\star)^n \subseteq s\left((J^\star)^n\right)^\star=
s\left(J^n\right)^\star\subseteq I^\star \, \Rightarrow \\
&  s((J^\star)^n \cap D) \subseteq s(J^\star)^n \cap D
\subset I^\star \cap D = I \subseteq P \, \Rightarrow \\ 
& (J^\star \cap D)^n \subseteq (J^\star)^n \cap D \subseteq P \, \Rightarrow \\ 
&  J^\star \cap D \subseteq P\,.
\end{array}
$$
Since $\star = \star_{f}$, then $P^\star = \cup \{ J^\star \mid \, J \in
{\boldsymbol {f}}(D)\,, \; J \subseteq P \}$ and so $P^\star \cap D = \cup
\{J^\star \cap D \mid \, J \in {\boldsymbol {f}}(D)\,, \; J \subseteq P \}
\subseteq P$; thus $P^\star \cap D = P$.  \hfill $\Box$

\vskip 5pt \hskip .5cm We denote by ${\cal M}(\star_{f}) \,$ \it the set
of all the quasi--$\star_{f}$--maximals of $D$\rm , which is nonempty if
and only if $\star_f \neq e$, and we associate to the semistar operation
$\star$ on $D$ a new semistar operation $ \tilde{\star}$ on $D$, which
is of finite type and spectral, defined as follows $ \tilde{\star} :=
\star_{{\cal M}(\star_{f}) }$\ (explicitly, $E^{\tilde{\star}} := \cap \{ED_Q
\,\mid\,\, Q \in {\cal M}(\star_f) \}\,$,\, for each $\,E \in
{\overline{\boldsymbol{F}}}(D)$\ ).\  Note that $\tilde{\star} \leq
\star_{f }$ \cite[Corollary 2.7]{FL3}.

\vskip 5pt \hskip .5cm We conclude this section by recalling the definition and the main
properties of the semistar Nagata rings.

\hskip .5cm  Let $D$ be an integral domain with field of quotients $K$ and $\star$ a
semistar operation on $D$.  Let $X$ be an indeterminate over $K$, for each $f\in
D[X]$, we denote by $\boldsymbol{c}(f)$ the content of $f$.  Let $N_{D}(\star):=\{h\in
D[X]\mid\, h\neq 0\, \,\hbox{and}\,\, \boldsymbol{c}(h)^\star=D^\star\}$.  Then
$N_{D}(\star)=D[X]\setminus\cup\{Q[X] \mid \, Q\in {\cal M}(\star_f)\}$ is a
saturated multiplicative system of $D[X]$.  The ring of fractions:
$$\hbox{Na}(D, \star):= D[X]_{N_{D}(\star)}$$
is called \it the Nagata ring of $D$ with respect to the semistar operation
$\star$ \rm (cf. \cite{FL3}). \rm

\vskip -3pt  Obviously, Na$(D, \star)=\hbox{Na}(D, \star_f)$ and if $\star=d$,
where $d$ is the identity (semi)star operation of $D$, then Na$(D, d)$
coincides with the ``classical'' Nagata ring $D(X)$ $:= \{f/g \mid f,
g\in D[X]\,,\; \boldsymbol{c}(g)=D\}$ of $D$.

%LEMMA 2.4
\begin{leem} \label{lm:2.4} \rm \cite [Corollary 2.7, Proposition 3.1
and 3.4, Corollary 3.5]{FL3}.  \sl Let $D$ be an integral domain with
quotient field $K$ and let $\star$ be a semistar operation on $D$. 
Then, for each $E\in {\overline{\boldsymbol{F}}}(D)$, we have: \balf

\rm \bf \item \sl 
$E^{\star_f}=\cap\{E^{\star_f}D_Q\mid\, Q\in {\cal M}(\star_f)\}$.

\rm \bf \item \sl 
$E^{\tilde{\star}}=\cap\{ED_Q\mid\, Q\in {\cal M}(\star_f)\}$.

\rm \bf \item \sl ${\cal M}(\star_f)={\cal M}({\tilde{\star}})$.

\rm \bf \item \sl  $\hbox{\rm Na}(D, \star)=\cap\{D[X]_{Q[X]}\mid\, Q\in
{\cal M}(\star_f)\}=\cap\{D_Q(X)\mid\,Q\in {\cal M}(\star_f)\}$.

\rm \bf \item \sl $\mbox{\rm Max}(\hbox{\rm Na}(D,
\star))=\{Q[X]_{N_{D}(\star)}\mid\, Q\in {\cal M}(\star_f)\}=
\{QD_Q(X) \cap \Na(D, \star) \mid\,Q\in {\cal M}(\star_f)\}$.

\rm \bf \item \sl 
$E\hbox{\rm Na}(D, \star)=\cap \{ED_Q(X)\mid\, Q\in {\cal M}(\star_f)\}$.

\rm \bf \item \sl 
$E\hbox{\rm Na}(D, \star)\cap K=\cap \{ED_Q\mid\, Q\in {\cal M}(\star_f)\}$.

\rm \bf \item \sl $E^{\tilde{\star}}=E\hbox{\rm Na}(D, \star)\cap K$.

\rm \bf \item \sl $\hbox{\rm Na}(D,
\star)=\hbox{\rm Na}(D, {\tilde{\star}})=\hbox{\rm Na}(D^{\tilde{\star}},
\dot{\tilde{\star}})$.  \hfill $\Box$

\ealf
\end{leem}

\hskip 0.5cm An easy consequence of the previous result (in particular, Lemma
\ref{lm:2.4} (e)) is the following:

%COROLLARY 2.5
\begin{coor} \label{cor:2.5} \sl Let $D$ be an integral domain and let
$\star$ be a semistar o\-pe\-ration on $D$.  For each
prime ideal $P$ of $D$ such that $P^\star \neq D^\star$, $\Na(D,
\star)_{P\Na(D, \star)}= D_{P}(X)$.
    
    \vskip -3pt \hfill $\Box$ \end{coor} \rm

%%%%%%%%%%%%%%%%%%%%%%%%%%%
    %%%%%%%%%%%%            SECTION 3
%%%%%%%%%%%%%%%%%%%%%%%%%%%

 %\begin{center}
    \section{Semistar linkedness}\label{sec:3}
    % \end{center}
    \vskip -9pt \hskip 0.5cm Let $D$ be an integral domain and $T$ be an
    overring of $D$.  Let $\star$ (respectively, $\star'$) be a semistar
    operation on $D$ (respectively, on $T$).

\vskip -3pt We say that $T$ is {\it $(\star, \star')$--linked to $D$} if:
$$
F^\star = D^\star \, \Rightarrow \, (FT)^{\star'}= T^ {\star'},$$
 for each nonzero finitely generated integral ideal $F$ of $D$.
 
  \vskip -3pt  It is straightforward that $T$ is $(\star, \star')$--linked to $D$
 if and only if $T$ is $(\star_{f}, \star'_{f})$--linked to $D$.
 
 \vskip -3pt  Obviously, $T$ is $(d_{D}, \star')$--linked to $D$, for each
  semistar operation $\star'$ on $T$ and $T$ is $(\star, e_{T})$--linked to
  $D$, for each semistar operation $\star$ on $D$; in particular, when $T
  $ coincides with the field of quotients $K$ of $D$, then there exists a
  unique (trivial) semistar operation $e_{T}= d_{T}$ on $T$, hence $T$ is
  $(\star, \star')$--linked to $D$, for each semistar operation $\star$ on
  $D$ and for each semistar operation $\star'$ on $T$.

 \vskip -3pt \hskip 0.5cm  We say that $T$ is {\it $t$--linked to } $(D, \star)$ if\ $T$ is $(\star,
  t_{T})$--linked.  In particular, the classical notion \it ``\ $T$ is
  $t$--linked to $D$'' \rm \cite{DHLZ1} coincides with the notion ``\ $T$ is
  $t$--linked to $(D, t_{D})$'' (i.e. $T$ is  $(t_{D}, t_{T})$--linked to $D$).

\smallskip

\hskip 0.5cm   In the following result we collect some
of the basic properties of the semistar linkedness.  

%LEMMA 3.1 
\begin{leem} \label{lm:3.1}  \sl Let $S, T$ be two overrings of an integral domain $D$,
with $D \subseteq T \subseteq S$. 
\balf
\rm\bf \item \sl Let $D = T$ and $\star', \star''$ be two semistar
operations on $T$.  If  $\star'_{f} \leq \star''_{f}$, then $T $ is $(\star',
\star'')$--linked to $T$.

\rm\bf \item \sl Let  $\star$ (respectively, $\star'$,
$\star''$) be a semistar operation on $D$ (respectively, $T$, $S$). 
Assume that $S$ is $(\star', \star'')$--linked to $T$ and that $T$ is
$(\star, \star')$--linked to $D$, then $S$ is $(\star, \star'')$--linked
to $D$.

\rm\bf \item \sl Let $\star$ (respectively, $\star'$, $\star''$) be a semistar operation
on $D$ (respectively, two semistar operations on $T$).  Assume that
$\star'_{f} \leq \star''_{f}$.  Then $T$ is $(\star, \star')$--linked to
$D$ implies that $T$ is $(\star, \star'')$--linked to $D$.

\rm\bf \item \sl  If $\star'$ is a (semi)star operation on $T$ (i.e. if\, 
$T^{\star'} =T$) and if\, $T$ is $(\star, \star')$--linked to
$D$ then $T$ is $t$--linked to $(D, \star)$.

\rm\bf \item \sl  Let $\star$ be a semistar operation on $D$ then $T$ is
$(\star, \dot{\star}^{T})$--linked to $D$.  In particular,
$D^{{\star}}$ is $({\star}, \dot{{\star}})$--linked to $D$.

\rm\bf \item \sl  If $\star'$ is a semistar operation on $T$ such that 
$\dot{\star}^{T}\le \star'$, then $T$ is $(\star, \star')$--linked to 
$D$. In particular,  we deduce that:

\smallskip

 \hskip 1cm ${\dot{\left(t_{D}\right)}}^{T}\leq {t}_{T} \, \Rightarrow \,
 T \,\, \hbox{is} \,\, t\hbox{--linked to} \, D \,;$ \;\; and more
 generally,\par

  \hskip 1cm $\left(\dot{\star}^{T}\right)_{f}\leq {t}_{T} \, \Rightarrow \, T
  \,\, \hbox{is} \,\, t\hbox{--linked to} \, (D, \star) \,$.

\rm\bf \item \sl Let $\star'$ be a semistar operation on $T$, then $T$ is
$(\hbox{\d{$\star$}}'_{D}, \star')$--linked to $D$.

\rm\bf \item \sl Let 
$\star_1$ and $\star_2$ be two semistar operations on $D$ and let $\star'$
be a semistar operation on $T$.
If $(\star_1)_{f}\leq (\star_2)_{f}$ and if  $T$ is $(\star_2,
\star')$--linked to $D$,  then $T$ is $(\star_1, \star')$--linked to $D$.

\rm\bf \item \sl  Let 
$\star$ (respectively, $\star'$) be a semistar operations on $D$
(respectively, $T$). If $\star \leq\hbox{\d{$\star$}}'_{D}$, then $T$ is
$(\star, \star')$--linked to $D$.

\vskip -2pt Also, we have:\par \hskip 1cm $t_{D} \leq {\hbox{\d{($t_{T}$}})}_{D}
\, \Rightarrow \, T \,\, \hbox{is} \, \, t\hbox{--linked to} \, D
\,;$\;\; and, more generally,

 \hskip 1cm  $\star_{f} \leq {\hbox{\d{($t_{T}$}})}_{D}   \, \Rightarrow \, T \,\, 
\hbox{is} \, \, t\hbox{--linked to} \, (D, \star)  \,$.

\bf \item \sl Let  $\star$ (respectively, $\star'$,
$\star''$) be a semistar operation on $D$ (respectively, $T$, $S$). 
Assume that $S$ is $(\star, \star'')$--linked to $D$ and that each 
quasi--$\star'_{f}$--maximal ideal
of $T$ is the contraction of a quasi--$\star''_{f}$--maximal ideal of 
$S$, then $T$ is $(\star, \star')$--linked to $D$.  

\noindent \vskip -5pt In particular (Lemma \ref{lm:2.3} (c)), if we take
$S:=T^{\star'}$ and $\star'':={\dot{\star'}}$ (note that
${\dot{\star'}}$ is a (semi)star operation on $T^{\star'}$), then $T$ is
$(\star, \star')$--linked to $D$ if and only if $T^{\star'}$ is $(\star,
{\dot{\star'}} )$--linked to $D$.

\bf \item \sl Let $\{ T_{\lambda}\,\mid\,\, \lambda \in \Lambda \}$ be a
family of overrings of $D$ and let $\ast_{\lambda}$ be a semistar operation
defined on $T_{\lambda}$, for $\lambda \in \Lambda $. 
Set $T:= \cap \{ T_{\lambda}\,\mid\,\, \lambda \in \Lambda \}$ and let
$\ast_{\Lambda}$ be the semistar operation on $T$ induced by the
family $\{ T_{\lambda}\,\mid\,\, \lambda \in \Lambda \}$ (i.e. for each 
$E\in \boldsymbol{{\overline{F}}}(T)$, $E^{\ast_{\Lambda}} := \cap \{
(ET_\lambda)^{\ast_{\lambda}} \mid \lambda \in \Lambda \}$).
If $T_{\lambda}$ is $(\star,
\ast_{\lambda})$--linked to $D$, for each $\lambda \in \Lambda$, then
$T$ is $(\star, \ast_{\Lambda})$--linked to $D$.

\ealf 
\end{leem} \rm

\bf Proof.  \rm Straightforward.  \hfill $\Box$
\medskip

\hskip .5cm Let $T, S$ be two overrings of an integral domain $D$, with
$D \subseteq T \subseteq S$ and let $\star$ (respectively, $\star'$,
$\star''$) be a semistar operation on $D$ (respectively, $T$, $S$). 
Assume that $S$ is $(\star, \star'')$--linked to $D$.  When is $S$
$(\star', \star'')$--linked to $T$ \hskip -1.5pt ?\ A partial answer to this
question will be given in  Remark \ref{rk:3.13}.

\medskip

%PROPOSITION 3.2
\begin{prro} \label{pr:3.2} \sl
Let $D$ be an integral domain and 
$T$ be an overring of $D$. Let $\star$ (respectively, $\star'$) be a 
semistar operation on $D$ (respectively,
on $T$).  The following are equivalent:\par
\brom 
\bf \item \sl $T$ is $(\star, \star')$--linked to $D$;\par

\bf \item \sl for each nonzero ideal $I$ of $D$, $I^{\star_{f}} = D^\star \,
\Rightarrow \, (IT)^{{\star'_f}} = T^{\star'}$;\par

\bf \item \sl for each quasi--$\star'_{f}$--ideal $J$ of $T$, with $J\neq 
T$,  $(J\cap D)^{\star_{f}} \neq D^\star$;\par

\bf \item \sl for each quasi--$\star'_{f}$--prime ideal $Q$ of $T$,
$(Q\cap D)^{\star_{f}} \neq D^\star $;\par

\bf \item \sl for each quasi--$\star'_{f}$--maximal ideal $N$ of $T$,
$(N\cap D)^{\star_{f}} \neq D^\star $.

\erom \end{prro} \rm

\bf Proof.  \rm  (i) $\Rightarrow$ (ii).  Since $D^\star= I^{\star_{f}} =
\cup\{F^\star \,| \; F \subseteq I\,, \; F \in \boldsymbol{f}(D)\}$,
then $D^\star =F^\star$, for some $F \subseteq I\,, \; F \in
\boldsymbol{f}(D)$.  Therefore, we conclude $T^ {\star'}=
(FT)^{\star'}\subseteq (IT)^{{\star'_f}} \subseteq T^ {\star'}$.
    
     (ii) $\Rightarrow$ (iii).  Assume that, for some
     proper quasi--$\star'_{f}$--ideal $J$ of $T$, the ideal $I:=J\cap D$ is
     such that $I^{\star_{f}} = D^\star$.  By assumption, we have $T^{\star'}
     = (IT)^{{\star'_f}} = ((J\cap D)T)^{{\star'_f}}\subseteq
     J^{{\star'_f}} \subseteq T^{\star'}$, i.e. $J^{{\star'_f}} =
     T^{\star'}$.  This fact contradicts the hypothesis that $J$ is a
     quasi--$\star'_{f}$--ideal of $T$, with $J \neq T$.

     (iii) $\Rightarrow$ (iv) $\Rightarrow$ (v) are obvious.
     
     (v) $\Rightarrow$ (i).  Assume that, for some $F \in \boldsymbol{f}(D)$,
     with $F \subseteq D$,  we have $F^\star = D^\star $ \, and \, $
     (FT)^{\star'_{f}}\subsetneq T^ {\star'}$.  Let $N$ be a
     quasi--$\star'_{f}$--maximal ideal of $T$ containing
     $(FT)^{\star'_{f}}\cap T$.  By hypothesis, we have $(N\cap
     D)^{\star_{f}} \neq D^\star $.  On the other hand, $F^{\star}\subseteq
     ((FT)^{\star'_{f}}\cap D)^{\star_{f}} \subseteq (N\cap D)^{\star_{f}}$
     and this contradicts the choice of $F$.  \hfill $\Box$

      %REMARK 3.3
\begin{reem} \label{rk:3.3}\rm \bf (a) \rm  It follows from Lemma \ref{lm:3.1} (b), (e) and (j)
     that, if $T^{\star'}$ is a $(\dot{\star}, \dot{\star'})$--linked \sl
     overring \rm of $D^{\star}$, then $T$ is $(\star, \star')$--linked to
     $D$.  What about the converse \hskip -2pt ?\ More precisely, since it is not true in
     general that $T^{\star'}$ is an overring of $D^{\star}$, for ``the
     converse'' we mean the following statement: 
     Assume $T^{\star'}$ is an overring of $D^{\star}$ 
     and that $T$ is $(\star, \star')$--linked to
     $D$.  Is it true that $T^{\star'}$ is $(\dot{\star},
     \dot{\star'})$--linked to $D^{\star}$ \hskip -1.5pt ?\   The answer to this question is
     negative, as the following example shows.
     
     \vskip -3pt   \hskip .5cm Let $K$ be a field and $X, Y$ be two indeterminates over $K$. Let 
$R:=K[X, Y]$ and $M:=(X, Y)$. Set $D:=K[X, XY]$ and $T:=R_M$. Let 
$\star:= 
{\hbox{\d{($t_{R}$}})}_{D}$ and $\star':=d_T$. Then:\par
\bara

\bf \item \rm $T$ is $(\star, \star')$--linked to $D$.

\bf \item \rm $D^{\star} \subseteq T^{\star'}$, but $T^{\star'}$ is not
$(\dot{\star}, \dot{\star'})$--linked to $D^{\star}$.  \eara

 \vskip -3pt    Clearly $D^{\star}=R\subseteq T^{\star'}=T$.

\vskip -3pt \bf (1) \rm Set $M':=MR_M$, then $M'$ is the unique 
($\star'$--)maximal ideal of $T$.  We have $M'\cap D=M\cap D\subseteq XR$. Therefore, $(M'\cap 
D)^{\star}\subseteq (XR)^{\star}=XR \subsetneq R= D^{\star}$.

\vskip -3pt \bf (2) \rm Note that $\dot{\star} =t_R$ (Example
\ref{ex:2.1} (6)) and  $\dot{\star}' = \star'=d_T$.  Moreover, for the
maximal ideal $M'$ of $T^{\star'} =T$, we have $(M'\cap
D^{\star})^{\dot{\star}}=(M'\cap R)^{t_R}=M^{t_R}=R=D^{\star}$. 
Therefore, $T^{\star'}$ is not $(\dot{\star}, \dot{\star'})$--linked to
$D^{\star}$ (Proposition \ref{pr:3.2} (v)).

     \vskip -3pt   \hskip .5cm  A related question to the previous one will be examined in Theorem \ref{th:3.7}.
     
     \bf (b) \rm If $T$ is $(\star, \star')$--linked to $D$, 
then, for each quasi--$\star'_{f}$--prime ideal $Q$ of $T$, there exists 
a quasi--$\star_{f}$--prime ideal $P$ such that $D_{P}
\subseteq T_{D\setminus P} \subseteq T_{Q}$.  (Since $(Q\cap D)^{\star_{f}} \neq
D^{\star}$, take a quasi--$\star_{f}$--prime ideal $P$ of $D$ such that
$Q\cap D \subseteq P$, and so $(D\setminus P) \subseteq (T\setminus
Q)$.)  Therefore, if $T$ is $(\star, \star')$--linked to $D$, then 
$D^{\widetilde{\star}} \subseteq T^{\widetilde{\star'}}$.

  \end{reem}

      %EXAMPLE 3.4
\begin{exxe} \label{ex:3.4} \rm \bf (1) \rm  
    Let $D$ be an integral domain and $T$ be an overring of $D$. 
 Let $\star$ be a semistar operation on $D$ and let $P$ be a
 quasi--$\star_{f}$--prime ideal of $D$.  Then, $T_{D \setminus P}$ is
 ${(\star, \ast)}$--linked to $D$, for \sl each \rm  semistar operation
 $\ast$ on $T_{D \setminus P}$ (equivalently, $T_{D \setminus P}$ is
 ${(\star, d_{T, P})}$--linked to $D$, where $d_{T, P}$ is the identity 
 (semi)star operation on $T_{D \setminus P}$). 
 
 \vskip -3pt As a matter of fact, for each prime
 ideal $N$, in particular, for each quasi--$\ast_{f}$--prime ideal, of
 $T_{D \setminus P}$, $N\cap T$ is a prime ideal of $T$ such that $N\cap
 D \subseteq P = P^{\star_{f}} \cap D$.  Hence $(N\cap D)^{\star_{f}}
 \neq D^\star$.

\bf (2) \rm Given a semistar operation $\star$ on an integral domain
$D$, recall that on $D$ we can introduce a new semistar operation of
finite type, denoted by $[\star]$, called
\it the semistar integral closure of $\star$,\rm $\,$ by setting: $$
F^{[\star]} := \cup \{((H^\star:H^\star)F)^{\star_f} \,\;|\;\, H \in
\boldsymbol{f}(D) \} \,,\;\;
\mbox{for each } \, F \in \boldsymbol{f}(D) \,, $$
(and thus in general:
$$E^{[\star]} := \cup \{F^{[\star]} \,\;|\;\, F \in
\boldsymbol{f}(D),\; F \subseteq E\}\,,\; \; \mbox{for each } \, E \in
\boldsymbol{\overline{F}}(D) \,).$$
It is known that $\, \star_{f} \leq [\star]\,$, hence $\, D^{\star}
\subseteq D^{[\star]}\,$, and that $D^{[\star]}$ is integrally closed. 
Therefore, it is obvious that if $\, D^{\star} = D^{[\star]}\, $ then
$D^{\star}$ is integrally closed.  The converse is false, even when $\,
\star\,$ is a (semi)star operation on $D$.  However, it is known that if
$\star_{f}$ is stable, then $D^{\star}$ is integrally closed if and only
if $\, D^{\star} = D^{[\star]}\,,$ \ (cf. 
\cite[Proposition 34]{Okabe/Matsuda:1994}, \cite[Proposition 4.3 and
Proposition 4.5]{FL1}, \cite [Example 2.1 (c)]{FJS}, \cite{Koch:1997}).

\vskip -4pt From Lemma \ref{lm:3.1} (e), (a) and (b), we have that
$D^{[\star]}$ is $(\star, \dot{[\star]})$--linked to $D$.

\end{exxe}
\medskip

 \hskip .5cm Assume that $T:= \cup \{ T_{\lambda} \mid \, \lambda \in \Lambda
 \}$ is the direct union of a given direct family of overrings $\{
 T_{\lambda}\;|\;\, \lambda \in \Lambda \}$ of an integral domain $D$
 with field of quotients $K$ (where $\Lambda$ is a directly ordered set
 by setting $\lambda' \leq \lambda''$ if $T_{\lambda'}\subseteq T_{\lambda''}$).  Let
 $\ast_{\lambda}$ be a semistar operation defined on the overring
 $T_{\lambda}$ of $D$, for each $\lambda \in \Lambda $.  We say that the
 family $\{\ast_{\lambda}\;|\;\, \lambda \in \Lambda \}$ is \it a direct
 family of semistar operations \rm (or, simply, that $\{(T_{\lambda},
 \ast_{\lambda})\mid \, \lambda \in \Lambda \}$ is a \it direct family\rm
 ), if $\lambda_{2}$ follows $\lambda_{1}$ inside $\Lambda$  and if $H\in
 \boldsymbol{f}(T_{\lambda_{1}})$, then $ H^{\ast_{\lambda_{1}}}
 \subseteq (HT_{\lambda_{2}})^{\ast_{\lambda_{2}}}$.
 
 \vskip -3pt For each $\lambda \in \Lambda$, let $E_{\lambda}$ be a
 $T_{\lambda}$-submodule of $K$.  We say that $E=\cup\{E_{\lambda}\mid\,
 \lambda\in \Lambda\}$ is \it a direct union\rm , if for each pair
 $\alpha, \beta \in \Lambda$, and for each $\gamma \in \Lambda$ such that
 $T_{\alpha} \subseteq T_{\gamma}$ and $T_{\beta} \subseteq T_{\gamma}$
 then $E_{\alpha}T_{\gamma} \subseteq E_{\gamma}$ and
 $E_{\beta}T_{\gamma} \subseteq E_{\gamma}$.
 
 \medskip
 
 \hskip .5cm The following result generalizes \cite[Proposition 2.2
 (a)]{DHLZ1}.

     %LEMMA 3.5
\begin{leem} \label{lm:3.5} \sl Let $\star$ be a semistar operation on
an integral domain $D$. Given a direct family $\{ (T_{\lambda}, \ast_{\lambda})
\;|\;\, \lambda \in \Lambda \}$, as above.  For each $E \in
\boldsymbol{\overline{F}}(T)$, set:
$$
E^{\ast^{\Lambda}} := \cup
\{E^{(\ast_{\lambda})_{f}} \mid \, \lambda \in \Lambda\} \,.
$$

\bara

\rm \bf \item\sl $\ast^{\Lambda}$ is a semistar operation of finite type
on $T$.

\rm \bf \item\sl If $T_{\lambda}$ is $(\star, \ast_{\lambda})$--linked to $D$,
for each $\lambda \in \Lambda $, then $T$ is $(\star,
\ast^{\Lambda})$--linked to $D$.

\rm \bf \item\sl If $T_{\lambda}$ is $(\star,
t_{T_{\lambda}})$--linked to $D$, for each $\lambda \in \Lambda $, then
$T$ is $(\star, t_{T})$--linked to $D$.

\eara

\end{leem}

\bf Proof.  \rm \bf (1) \rm The properties $(\star_1)$ and $(\star_2)$ 
are straightforward.
 Before proving 
$(\star_3)$, we show the following:

\bf Claim. \sl  If $E=\cup\{E_{\lambda}\mid\, 
\lambda\in \Lambda\}\in 
\overline{\boldsymbol{F}}(T)$ is  a direct  union, where $E_{\lambda}$ is a  
$T_{\lambda}$-submodule of $K$, then: 
$$E^{\ast^{\Lambda}}=\cup\{E_{\lambda}^{({\ast}_{\lambda})_f}\mid\, 
\lambda\in \Lambda\}\,.$$

\rm Given $\alpha \in \Lambda$, we have $E=\cup\{E_{\beta}T_{\alpha}\mid\, 
\beta\in\Lambda\}$ is a direct union of  $T_{\alpha}$-submodules. Since $({\ast}_{\alpha})_f$ is of finite 
type and $E \in 
\overline{\boldsymbol{F}}(T_{\alpha}) \ (\supseteq 
\overline{\boldsymbol{F}}(T)) $, then  
$E^{(\ast_{\alpha})_f}=\cup\{({E_{\beta}T_{\alpha}})^{({\ast}_{\alpha})_f}
\mid\, \beta\in \Lambda\}$. Let $\beta\in \Lambda$, then there exists  $\gamma\in\Lambda$ such 
that $T_{\alpha}\subseteq T_{\gamma}$ and $T_{\beta}\subseteq T_{\gamma}$ 
and, $E_{\alpha}\subseteq E_{\gamma}$ and $E_{\beta}\subseteq E_{\gamma}$. 
Hence $({E_{\beta}T_{\alpha}})^{({\ast}_{\alpha})_f}\subseteq 
{E_{\gamma}}^{({\ast}_{\alpha})_f}\subseteq 
{E_{\gamma}}^{({\ast}_{\gamma})_f}$ (the second inclusion follows from 
the fact that $\{ (T_{\lambda}, \ast_{\lambda}) \,\mid\, \lambda\in 
\Lambda \}$ is direct). So 
$E^{(\ast_{\alpha})_f}\subseteq\cup\{E_{\lambda}^{({\ast}_{\lambda})_f}\mid\, 
\lambda\in \Lambda\}$, and hence 
$E^{\ast^{\Lambda}}\subseteq\cup\{E_{\lambda}^{({\ast}_{\lambda})_f}\mid\, 
\lambda\in \Lambda\}$. The other inclusion is trivial.

\vskip -3pt Now we prove $(\star_3)$.  Clearly, for each $E \in
\overline{\boldsymbol{F}}(T)$, $E\subseteq E^{\ast^{\Lambda}}$.  On the
other hand, we have $E^{\ast^{\Lambda}} = \cup
\{E^{(\ast_{\lambda})_f}\, \mid \, \lambda \in \Lambda \}$ is a direct
union of $E^{(\ast_{\lambda})_f } \in
\overline{\boldsymbol{F}}(T_{\lambda})$ and so, by the Claim,
$(E^{\ast^{\Lambda}})^{\ast^{\Lambda}} = \cup
\{(E^{(\ast_{\lambda})_f})^{(\ast_{\lambda})_f}\,\mid\, \lambda
\in\Lambda\}=\cup \{E^{(\ast_{\lambda})_f}\,\mid\, \lambda
\in\Lambda\}=E^{\ast^{\Lambda}}$.

\vskip -3pt Finally, the fact that ${\ast}^{\Lambda}$ is of finite type is an immediate  
consequence of the definition. 

\bf (2) \rm Let $I$ be a nonzero finitely generated ideal of
$D$ such that $I^\star = D^\star$.  Then, by the Claim, 
$(IT)^{\ast^{\Lambda}} = \cup \{(IT_{\lambda})^{\ast_{\lambda}}\,\mid\,
\lambda \in \Lambda\}$.  Since $T_{\lambda}$ is $(\star,
\ast_{\lambda})$--linked to $D$, then
$(IT_{\lambda})^{\ast_{\lambda}}=T^{\ast_{\lambda}}_\lambda$, for each
$\lambda \in \Lambda$.  Hence, again by the Claim, 
$(IT)^{\ast^{\Lambda}}= \cup \{T_{\lambda}^{\ast_{\lambda}}\,\mid\,
\lambda \in \Lambda\}= T^{\ast^{\Lambda}}$.  \par

 \bf (3) \rm Let $I$ be a nonzero finitely generated ideal of $D$ such that
 $I^\star = D^\star$, then for each $\lambda$, $(IT_{\lambda})^t =
 T_{\lambda}$, i.e. $(IT_{\lambda})^{-1}= T_{\lambda}$.  Let $I :=
 (x_{1}, x_{2}, \ldots, x_{n})D$ and $z \in (IT)^{-1}$.  Then, for each
 $i$, $zx_{i}\in T_{\lambda_{i}}$, for some $\lambda_{i}\in \Lambda$ and
 so, for some $\lambda_{I}\in \Lambda$, $zI \subseteq T_{\lambda_{I}}$. 
 Hence, $z\in (IT_{\lambda_{I}})^{-1} =T_{\lambda_{I}}\subseteq T$. 
 Therefore, $(IT)^{-1} \subseteq T$ and so $(IT)^{-1} =T$.  \hfill $\Box$
 
 \medskip
 \hskip 0.5cm The following corollary generalizes \cite[Corollary
 2.3]{DHLZ1}.

   %COROLLARY 3.6
\begin{coor} \label{cor:3.6} \sl Let $\star$ be a semistar operation
on an integral domain $D$. Then $D^{[\star]}$ is $t$--linked to $(D, \star)$. 
If, moreover $(D:D^\star ) \neq (0)$, then the complete integral closure
$\widetilde{D}$ of $D$ is $t$--linked to $(D, \star)$; in particular,
the complete integral closure $\widetilde{D}$ of $D$ is always
$t$--linked to $D$.
    \end{coor}
    
    \bf Proof.  \rm The statement can be seen as an easy consequence of
    Example \ref{ex:3.4} (2) and of the fact that $\dot{[\star]} \le
    t_{D^{[\star]}}$ (Lemma \ref{lm:3.1} (c)).  We give here another proof
    based on the previous Lemma \ref{lm:3.5}, which also shows that the
    semistar operation $[\star]$ is issued from a semistar operation
    associated to a directed family of overrings and semistar operations.
    
    \vskip -3pt \hskip 0.5cm For each $E\in {\overline{\boldsymbol{F}}}(D)$,
    set $T_E :=(E^\star:E^\star)$.  Let $\ast_E$ denote the semistar operation
    $\dot{\star}^{T_E}$ on $T_E$.  Then $T_E$ is an overring of $D$, which is
    $(\star, \ast_E)$--linked to $D$ (Lemma \ref{lm:3.1} (e)).  Note that
    $\ast_E$ is a (semi)star operation on $T_E$ (since
    $(T_{E})^{\ast_{E}}=T_E$).  
    
    \vskip -3pt  We claim that $\{(T_F, \ast_F)\mid\, F\in
    {\boldsymbol{f}}(D)\}$ and $\{(T_E, \ast_E)\mid\, E\in
    {{\boldsymbol{F}}}(D)\}$ are direct families  (as in Lemma
    \ref{lm:3.5}).  To see this,
note that:
$$
(H_{1}^\star:H_{1}^\star) \subseteq
((H_{1}H_{2})^\star:(H_{1}H_{2})^\star) \supseteq 
(H_{2}^\star:H_{2}^\star),
$$
for all $H_{1}, H_{2}\in {\overline{\boldsymbol{F}}}(D)$.

\vskip -3pt Therefore, as in Lemma \ref{lm:3.5} (1),  $\{ (T_F, \ast_F)\mid\, F\in
    {\boldsymbol{f}}(D) \}$ (respectively, $\{ (T_E, \ast_E) \mid\, E\in
    {{\boldsymbol{F}}}(D) \}$) defines a (semi)star operation of finite type
    $\ast^{{\boldsymbol{f}}(D)}$ (respectively,
    $\ast^{{\boldsymbol{F}}(D)}$) on $D^{[\star]} = \cup \{ (F^\star:F^\star)
    \mid \, F \in {\boldsymbol{f}}(D) \}$ (respectively, on $D^{\langle
    \star \rangle} := \cup \{ (E^\star:E^\star) \mid \, E \in
    {\boldsymbol{F}}(D)\}$).

\vskip -3pt Note that  $D^{[\star]}$ is $(\star,
\ast^{{\boldsymbol{f}}(D)})$--linked to $D$ (Lemma \ref{lm:3.5} (2)) and 
that $\ast^{\boldsymbol{f}(D)}\leq t_{D^{[\star]}}$  (since  $\ast^{\boldsymbol{f}(D)}$ 
is a (semi)star operation of finite type on $D^{[\star]}$).  We conclude,
by Lemma \ref {lm:3.1} (c), that $D^{[\star]}$ is $t$--linked to $(D,
\star)$.

\vskip -3pt \hskip 0.5cm For the last statement, note that
$\widetilde{D}= \cup \{(E:E)\,\mid\, E \in{{{\boldsymbol{F}}}(D)} \}
\subseteq \cup \{(E^\star:E)\,\mid\, E \in {{{\boldsymbol{F}}}(D)} \} =
\cup \{(E^\star:E^\star)\,\mid\, E \in {{{\boldsymbol{F}}}(D)} \} =
D^{\langle \star \rangle} \subseteq \cup \{(H:H)\,\mid\, H \in
{{{\boldsymbol{F}}}(D^\star)} \} = \widetilde{D^\star}$.  If $(D :D^\star) \neq (0)$ then $\widetilde{D}=\widetilde{D^\star}= D^{\langle \star \rangle} =\cup \{T_E\,\mid\,
E \in {{\boldsymbol{F}}}(D) \}$.  Arguing as above, we have that
$\widetilde{D}$ is $(\star, \ast^{{{\boldsymbol{F}}}(D)})$--linked to
$D$, and $\ast^{\boldsymbol{F}(D)} \leq
t_{\widetilde{D}}$ (since $\ast^{\boldsymbol{F}(D)}$ 
is a (semi)star operation of finite type on $D^{\langle \star \rangle}=
\widetilde{D}$).  Again from Lemma \ref{lm:3.1} (c), we conclude that
$\widetilde{D}$ is $t$--linked to $(D, \star)$.  \hfill $\Box$

%REMARK 3.7
\begin{reem} \label{rk:3.7} \rm Let $\star$ be a semistar operation on an
integral domain $D$.

\vskip -2pt \bf (a) \rm Let $\ast^{\boldsymbol{f}(D)}$ be  the
(semi)star operation of finite type over $D^{[\star_{_{\!f}}]}\ 
(=D^{[\star]})$, associated to the semistar operation of finite type 
$\star_{f}$ and defined, in general for any semistar operation,  in
the proof of the previous corollary. Then:
$$
[\star] = \hbox{\d{$\left(\ast^{{\boldsymbol{f}}(D)}\right)$}}_{D}\,.
$$
As a matter of fact, first, note that in this case $T_H= (H^{\star_{_{\!f}}}:
H^{\star_{_{\!f}}}) = (H^\star : H^\star) $, for each $H \in
{{\boldsymbol{f}}}(D)$ and
 let now $\ast_H$ denote the semistar operation of finite type
    $\dot{\star}_{_{\!f}}^{T_H}$ on $T_H$.  For each $E \in
    {\overline{{\boldsymbol{F}}}}(D)$, we have:
 $$
ED^{[\star]} = E\left( \cup \{ (H^\star : H^\star) \mid \, H\in
{{\boldsymbol{f}}}(D) \} \right) = \cup \{ E(H^\star : H^\star) \mid \,
H\in {{\boldsymbol{f}}}(D) \}\,,
$$
 thus, using the Claim of the proof of Lemma \ref{lm:3.5}, we have:
		$$
\begin{array}{rl}
   (ED^{[\star]})^{\ast^{{\boldsymbol{f}}(D)}} &= \left( \cup \{ E(H^\star
   : H^\star) \mid \, H\in {{\boldsymbol{f}}}(D) \}
   \right)^{\ast^{{\boldsymbol{f}}(D)}} =\\
&= \cup \{ 
    \left(E\left( H^\star :  H^\star \right)\right)^{\star_{f}} 
    \mid \, H \in {{\boldsymbol{f}}}(D) \} \,.
    \end{array}
    $$
    In particular, $F^{[\star]}  = (FD^{[\star]})^{\ast^{{\boldsymbol{f}}(D)}}$,\ for each $F\in {{\boldsymbol{f}}}(D)$.
    
    As a consequence we have that, for each $E \in
    {\overline{\boldsymbol{F}}}(D)$:
    $$
    E^{[\star]} = \cup \{ 
    \left(E\left( H^\star : H^\star \right)\right)^{\star_{f}} \mid \, H \in
    {{\boldsymbol{f}}}(D) \}\,.
    $$
    
    \bf (b) \rm If we set:
    $$
    \langle \star \rangle :=
    \hbox{\d{$\left(\ast^{{\boldsymbol{F}}(D)}\right)$}}_{D}\,,
$$
then $
    \langle \star \rangle$ is a  semistar operation of finite type on
    $D$, with $D^{\langle \star \rangle}= \cup \{(E^\star:E^\star)\,\mid\, E
    \in {{{\boldsymbol{F}}}(D)} \}$.  Moreover,
$$
\star_{f}\leq [\star] \leq \langle \star \rangle \; \; \mbox{ and } \;\;
D^\star \subseteq D^{[\star]} \subseteq D^{\langle \star \rangle}
\subseteq \widetilde{D^{\star}}\,.
$$
\end{reem}
    
    %THEOREM 3.8
\begin{thee} \label{th:3.7} \sl Let $D$ be an integral domain with quotient field $K$
    and let $T$ be an overring of $D$.  Let $\star$ (respectively, $\star'$)
    be a semistar operation on $D$ (respectively, on $T$).  The following
    are equivalent:
    \brom
    
   \rm \bf \item \sl  $T$ is $(\star, \star')$--linked to $D$;
   
    \rm \bf \item \sl $\Na(D, \star) \subseteq \Na(T, \star')$;

     \rm \bf \item \sl ${\tilde{\star}}\le \hbox{\d{$(\widetilde{\star'})$}}_{D}$;

     \rm \bf \item \sl $T$ is an $({\tilde{\star}},
     {\widetilde{\star'}})$--linked overring of $D$;

    \rm \bf \item \sl  $T^{\widetilde{\star'}}$ is an $(\dot{\tilde{\star}},
    \dot{\widetilde{\star'}})$--linked overring of $D^{\tilde{\star}}$.
    
    \erom
    \end{thee}
    
     \bf Proof.  \rm  (i) $\Rightarrow$ (ii).  Let $g\in D[X]$
    such that $(\boldsymbol{c}_D(g))^\star=D^\star$.  Then, by the assumption,
    $(\boldsymbol{c}_T(g))^{\star'}=(\boldsymbol{c}_D(g)T)^{\star'}=T^{\star'}$.  Hence Na$(D,
    \star)\subseteq \hbox{Na}(T, \star')$.
    
    (ii)$ \Rightarrow $(iii).  Let
    $E\in {\overline{\boldsymbol{F}}}(D)$.  Then $E\hbox{Na}(D, \star)\subseteq
    E\hbox{Na}(T, \star')$.  Hence (Lemma \ref{lm:2.4} (h)) $E^{\tilde{\star}}=
    E\hbox{Na}(D, \star)\cap K\subseteq E\hbox{Na}(T, \star')\cap
    K=(ET)^{\widetilde{\star'}}$ and so we conclude that ${\tilde{\star}}\le
    \hbox{\d{$(\widetilde{\star'})$}}_{D}$\,.
    
    (iii) $ \Rightarrow$ (iv).  It
    follows from Lemma \ref{lm:3.1} (i).

    (iv) $\Rightarrow$ (ii) follows from (i) $\Rightarrow$ (ii) and from
    Lemma \ref{lm:2.4} (i).

    (ii) $\Rightarrow$ (i).  Let $G$ be a nonzero finitely generated integral
    ideal of $D$ such that $G^\star = D^\star$ and let $g\in D[X]$ be such
    that $\boldsymbol{c}_D(g) = G$.  From the fact that
    $(\boldsymbol{c}_D(g))^\star=D^\star$, we have that $g$ is a unit in Na$(D,
    \star)$ and so, by   assumption, $g$ is also a unit in Na$(T, \star')$. 
    This implies that
    $(\boldsymbol{c}_T(g))^{\star'}= (\boldsymbol{c}_D(g)T)^{\star'}=T^{\star'}$,
    i.e. $(GT)^{\star'} =T^{\star'}$.
    
    (ii) $\Leftrightarrow$ (v) is an easy consequence of (ii)
    $\Leftrightarrow$ (i) and of Lemma \ref{lm:2.4} (i).
   \hfill
    $\Box$
    \medskip

    \hskip .5cm The next result characterizes domains such that each overring
    is semistar linked and generalizes \cite[Theorem 2.6]{DHLZ1}.

\medskip

    %THEOREM 3.9
\begin{thee} \label{th:3.8} \sl Let $D$ be an integral domain and
$\star$ a semistar operation on $D$.  The following statements are equivalent:
\brom
 \rm \bf \item \sl 
   for each overring $T$ of $D$ and for each semistar operation
$\star'$ on $T$, $T$ is $(\star, \star')$--linked to $D$;

\rm \bf \item \sl 
each overring $T$ of $D$ is $(\star, d_{T})$--linked to $D$;

\rm \bf \item \sl each
overring $T$ of $D$ is $t$--linked to $(D, \star)$;

\rm \bf \item \sl for each
valuation overring $V$ of $D$ there exists a (semi)star operation
$\ast_{V}$ on $V$, such that $V$ is $(\star, \ast_{V})$--linked to
$D$;

\rm \bf \item \sl each maximal ideal of $D$ is a quasi--$\star_{f}$--maximal 
ideal;

\rm \bf \item \sl for each  proper ideal $I$ of $D$, $I^{\star_{f}}
\subsetneq D^\star$;

\rm \bf \item \sl for each proper finitely generated ideal $I$ of $D$,
$I^{\star} \subsetneq D^\star$;

\rm \bf \item \sl  for each proper $\star_f$--invertible ideal $I$ of $D$ (i.e.
$(II^{-1})^{\star_f} = D^\star$), $I^{\star_f} \subsetneq D^\star$
(hence, each proper $\star_f$--invertible ideal $I$ of $D$ is contained
in the proper quasi--$\star_{f}$--ideal $I^\star \cap D$ of $D$).  \erom
\end{thee}

 \bf Proof.  \rm   (i) $\Rightarrow$ (ii) is obvious. 
(ii) $\Rightarrow$ (iii) is a consequence of the fact that $d_T\le t_T$
and Lemma \ref{lm:3.1} (c).  (iii) $\Rightarrow$ (iv) is obvious, taking
$\ast_{V}= t_{V}$.

(iv) $\Rightarrow$ (v).  If $M$ is a maximal
ideal of $D$ such that $M \subsetneq M^{\star_{f}}=D^\star$ then, for some
nonzero finitely generated ideal $I \subseteq M$, we have $I^\star=
D^\star$.  Let $(V, N)$ be a valuation overring of $D$ such that $N \cap
D = M$.  Then $(IV)^{(\ast_{V})_{f}} = V^{\ast_{V}} = V$.  Since a
nonzero finitely generated ideal of a valuation domain is principal and
$\ast_{V}$ is a (semi)star operation on $V$, then $V = V^{\ast_{V}} =
(IV)^{(\ast_{V})_{f}} = IV$.  This is a contradiction, because
$IV\subseteq N \subsetneq V$.

(v) $\Rightarrow$ (vi) $\Rightarrow
$(vii) are obvious.

(vii) $\Rightarrow$ (viii).  If
$(II^{-1})^{\star_f} = D^\star$ and $II^{-1}\subsetneq D$ then, for some
nonzero finitely generated ideal $F \subseteq II^{-1}\subsetneq D$, we
have $F^\star = D^\star$ and this contradicts the assumption.  Since $I$ is
invertible then, in particular, $I$ is a finitely generated proper ideal
of D and so, by assumption and (vii), $I^\star \cap D$ is a proper
quasi--$\star_{f}$--ideal of $D$ containing $I$.  \par

(viii) $\Rightarrow$ (v). Assume that, for some maximal ideal $M$ of 
$D$, $M \subsetneq M^{\star_{f}} = D^\star$.  Then
$(MM^{-1})^{\star_{f}}= (MM^{-1})^\star = D^\star$, because $
D^{\star_{f}} \supseteq (MM^{-1})^{\star_{f}} =
(M^{\star_{f}}(M^{-1})^{\star_{f}})^{\star_{f}}
=(D^{\star_{f}}(M^{-1})^{\star_{f}})^{\star_{f}} = (M^{-1})^{\star_{f}}
\supseteq D^{\star_{f}}= D^\star$.  Hence, by assumption,
$M^{\star_f}\subsetneq D^\star$, but this contradicts the choice of $M$.

(v) $\Rightarrow$ (i).  Assume that, for some overring $T$
of $D$, for some semistar operation $\star'$ on $T$ and for some 
quasi--$\star'_{f}$--maximal ideal $N$ of $T$, we have $(N\cap
D)^{\star_{f}}=D^\star$ (Proposition \ref{pr:3.2} ((i) $\Leftrightarrow$ (v))). 
Note that, from the assumption, $N\cap D \subseteq M = M^{\star_{f}}
\cap D$, for some (quasi--$\star_{f}$--)maximal ideal $M$ of $D$, and so
we reach immediately a contradiction.  \hfill $\Box$

%REMARK 3.10
\begin{reem} \label {rk:3.10} \rm Note that  the proof of (vii)
$\Rightarrow$ (viii) (Theorem \ref{th:3.8}) shows that, in an integral
domain verifying the conditions of Theorem \ref{th:3.8}, each
$\star_f$--invertible ideal is invertible.
    \end{reem}

%EXAMPLE 3.11
\begin{exxe} \label {ex:3.9} \rm Let $\star$ be a semistar operation on
an integral domain $D$.  Assume that $D^\star$ is faithfully flat on
$D$ (for instance, assume that $\star$ is a (semi)star operation on $D$).
In this situation, every principal ideal of $D$ is a quasi--$\star$--ideal of
$D$.  If $\hbox{Spec}(D)$ is a tree (e.g., dim$(D)=1$ or $D$ is a
GD-domain, in particular, $D$ is a Pr\"ufer domain), then every
overring $T$ of $D$ is $t$--linked to $(D, \star)$.

\vskip -3 pt \hskip .5cm In order to apply Theorem \ref{th:3.8} ((v)
$\Rightarrow$ (iii)), we show that each maximal ideal $M$ of $D$ is a
quasi--$\star_{f}$--ideal of $D$.  For each nonzero $x \in M$, $xD$ is a
quasi--$\star_{f}$--ideal of $D$, hence a minimal prime ideal $P$ of
$xD$ is a quasi--$\star_{f}$--prime ideal of $D$ (cf.  Lemma
\ref{lm:2.3} (d)).  Since $\hbox{Spec}(D)$ is a tree, $M$ is a direct
union of a family $\{P_{\lambda}\}$ of quasi--$\star_{f}$--prime ideals
of $D$.  If $M^{\star_{f}} = D^\star$, then $1 \in M^{\star_{f}} =
(\cup_{\lambda} \{P_{\lambda}\})^{\star_{f}} = (\cup_{\lambda}
\{(P_{\lambda})^{\star_{f}}\})^{\star_{f}}$ thus, from the finiteness of
${\star_{f}}$, we deduce that $1 \in (P_{\lambda})^{\star_{f}}\cap D =
P_{\lambda}$, for some $\lambda$, and this is a contradiction.

\end{exxe}

\medskip
\hskip .5cm Our next goal is the study of a new semistar operation strictly
related to semistar linkedness.

\vskip -3 pt  Let $D$ be an integral domain,
$\star$ a semistar operation on $D$, and $T$ an overring of $D$.  We define
\it the semistar operation ${\ell}_{^{\star, T}}$\rm (or, simply,
$\ell$) \it on \rm $T$, in the following way:
$$
E^{{\ell}_{^{\star, T}}} := E^{\ell} := \cap\{ET_{D\setminus P}\, | 
\,\, P
\, \hbox{is a quasi--$\star_{f}$--prime ideal of } \, D \}\,,
$$
for each $E \in {\overline{\boldsymbol{F}}}(T) $.

\vskip -3 pt  Note that if $T = D$, then
${\ell}_{^{\star, D}} = \tilde{\star}$ (Lemma \ref{lm:2.4} (b)). 
Moreover, note that ${\ell}_{^{\star, T}}$ is the semistar operation on
$T$ induced, in the sense described in Example \ref{ex:2.1} (4), by the family
of overrings $\{T_{D\setminus P}\mid\, P\,\, \hbox{is a
quasi--$\star_{f}$--prime ideal of } \, D \}$ of $D$ (where
$T_{D\setminus P}$ is endowed with the identity $d_{T, P}$ (semi)star
operation) .

\hskip .5cm The following proposition collects some interesting
properties of the semistar operation ${\ell}_{^{\star, T}}$ .

    %PROPOSITION 3.12
\begin{prro} \label{pr:3.10}  \sl
    Let $D$ be an integral domain, $\star$ a
semistar operation on $D$, $T$ an overring of $D$ and $\star'$ a semistar
operation on $T$.
\bara \rm \bf \item \sl ${\ell}_{^{\star,T}}$ is a stable semistar
operation of $T$.

\rm \bf \item \sl   Assume that $T$ is
$(\star, \star')$--linked to $D$.  Then\ ${\ell}_{^{\star,T}}
\leq\widetilde{\star'} \ (\leq \star'_{f})$; in parti\-cu\-lar $T$ is
$({\ell}_{^{\star, T}}, \star')$--linked to $T$.

\rm \bf \item \sl   $T$ is $(\star,
{\ell}_{^{\star, T}})$--linked to $D$, for each semistar operation 
$\star$
on $D$;  in particular, $D$ is $(\star,
\tilde{\star})$--linked to $D$, for each semistar operation $\star$ 
on $D$.

\rm \bf \item \sl   ${\ell}_{^{\star, T}}$ is a semistar operation of finite type on 
$T$ and $\widetilde{{\ell}_{^{\star,T}}}={\ell}_{^{\star,T}}$. 

\rm \bf \item \sl    $\ell_{\star, T}$ is the unique 
minimal element in  set of semistar operations $\star'_f$, where 
$\star'$ is a semistar operation  on $T$ such that $T$ is $(\star, 
\star')$--linked to $D$. 

\rm \bf \item \sl   $T$ is $(\star, \star')$--linked to $D$ if and
only if $T$ is
$({\ell}_{^{\star, T}}, \star')$--linked to $T$ (and $T$ is
$(\star, {\ell}_{^{\star, T}})$--linked to $D$).

\rm \bf \item \sl   $T$ is $(\star, \star')$--linked to $D$ if and
only if\, ${\ell}_{^{\star, T}} \leq \star'_f$.

\eara

\end{prro}

\bf Proof.  \rm \bf (1) \rm  This is a straightforward consequence of the
fact that $T_{D\setminus P}$ is flat over $T$, for each prime ideal $P$
of $D$.

\bf (2) \rm For each quasi--$\star'_{f}$--prime ideal $Q$ of $T$, there
exists a quasi--$\star_{f}$--prime ideal $P$ of $D$, such that
$T_{D\setminus P} \subseteq T_{Q}$ (Remark \ref{rk:3.3} (b)), and so
also $ET_{D\setminus P} \subseteq ET_{Q}$, for each $E \in
{\overline{\boldsymbol{F}}}(T)$; from this we deduce that
${\ell}_{^{\star, T}} \leq \widetilde{\star'}$.  The last statement
follows from Lemma \ref{lm:3.1} (a). 

\bf (3) \rm If $I^{\star_{f}}= D^\star$, then
$I\not\subseteq P$, i.e. $ID_{P}=D_{P}$, and this implies that $IT_{D
\setminus P}=T_{D\setminus P}$, for each quasi--$\star_{f}$--prime ideal
$P$ of $D$.  Therefore $(IT)^{{\ell}_{^{\star, T}} } =
T^{{\ell}_{^{\star, T}}}$.  

\bf (4) \rm From (3), we have that $T$ is $(\star, {\ell}_{^{\star,
T}})$--linked to $D$.  From (2) (for $\star' = {\ell}_{^{\star,T}}$), we deduce that
$({\ell}_{^{\star,T}})_{f}\le{\ell}_{^{\star,T}}
\leq\widetilde{{\ell}_{^{\star,T}}} \leq ({\ell}_{^{\star,T}})_{f}$.\par
\bf (5) \rm follows from (2) and (3).

\bf (6) \rm It is a direct
consequence of (2), (3) and Lemma \ref{lm:3.1} (b). 

\bf (7) \rm is equivalent to (6), by (2) and Lemma \ref{lm:3.1} (a).  \hfill $\Box$

\bigskip

    %REMARK 3.13
\begin{reem} \label{rk:3.13} \rm Let $T, S$ be two overrings of an integral
domain $D$, with $D \subseteq T \subseteq S$ and let $\star$
(respectively, $\star'$, $\star''$) be a semistar operation on $D$
(respectively, $T$, $S$).  Assume that $S$ is $(\star, \star'')$--linked
to $D$.  If $T$ is $(\star', {\ell}_{^{\star,T}})$--linked to $T$ (e.g.
if $\star'_{f} \leq {\ell}_{^{\star,T}}$), then $S$ is $(\star',
\star'')$--linked to $T$.  As a matter of fact, let $Q$ be a
quasi-$\star''_f$-prime ideal of $S$, then $(Q\cap D)^{\star_{f}} \neq
D^\star$, and hence, by definition of  $T^{{\ell}_{^{\star, T}}}$, $(Q\cap T)^{{\ell}_{^{\star, T}}} \neq
T^{{\ell}_{^{\star, T}}}$.  So $(Q\cap T)^{\star'_{f}} \neq T^{\star'}$.
 
\vskip -3pt \hskip .5cm In general, for \sl any \rm nontrivial semistar
operation $\star''$ on $S$, we can construct a nontrivial semistar
operation $\star'$ on $T$ such that $S$ is not $(\star',
\star'')$-linked to $T$: Let $Q$ be a quasi--$\star''_f$--prime ideal of
$S$, and let $0\not=q\in Q\cap T$.  Let $T_{q}$ be the ring of fractions
of $T$ with respect to its multiplicative set $\{q^n\mid\, n\ge 0\}$ and
let $\star':=\star_{\{T_{q}\}}$.  Then $S$ is not $(\star',
\star'')$-linked to $T$, since $(Q\cap T)^{\star'_{f}}=(Q\cap
T)T_q=T_q=T^{\star'}$.
    \end{reem}

\hskip .5cm  In \cite{DHLZ1}, the authors showed that the equality
$T^{{\ell}_{t_D,T}}=T$ characterizes $t$--linkedness of $T$ to $D$. 
The next goal is to investigate the analogous question in semistar
setting.\par

    %LEMMA 3.14
\begin{leem} \label{lm:3.11} \sl Let $D$ be an integral domain, $T$an
overring of $D$, $\star$  a semistar operation on $D$ and
$\star'$ a (semi)star operation on $T$.  If $T$ is $(\star,
\star')$--linked to $D$, then $T^{{\ell}_{\star,T}}=T$. 
\end{leem}
\bf Proof.  \rm Since $\star'$ is a (semi)star operation on $T$, then
$T =T^{\widetilde{\star'}}=T^{\star'}$.  Therefore, by Proposition \ref{pr:3.10} (2),
we have $T\subseteq T^{{\ell}_{\star,T}}\subseteq
T^{\widetilde{\star'}}=T$, and so $T^{{\ell}_{\star,T}}=T$.  \hfill
$\Box$

\medskip

\hskip .5cm However, ``a general converse'' of the previous lemma  fails to be
true as the following example shows.

   %EXAMPLE 3.15
\begin{exxe} \label{ex:3.12} \rm  Let $K$ be a field and $X, Y$ two
indeterminates over $K$.  Let $D:=K[X, Y]$ and $M:=(X, Y)$.  Set $T:=
D_M$.  Then $D\subset T$ is $t$--linked (since $D\subset  T$ is flat
\cite[Proposition 2.2 (c)]{DHLZ1}).  Hence $T^{{\ell}_{^{{t_D}, T}}}=T$,
by \cite[Proposition 2.13 (a)]{DHLZ1}.  On the other hand, we have
$MT\neq T$ and $M^{t_D}=D$.  Hence $T$ is not $(t_D, d_T)$--linked to
$D$.
\end{exxe}

\hskip .5cm  A generalization of \cite[Proposition 2.13 (a)]{DHLZ1} is
given next, by showing that the converse of Lemma \ref{lm:3.11} holds
when $\star'=t_T$.

   %PROPOSITION 3.16
\begin{prro} \label{pr:3.13} \sl  Let $\star$ be a semistar operation on 
the integral domain $D$ and $T$ an overring of $D$. Then 
 $T$ is $t$--linked to 
$(D, \star)$ if and only if $T^{{\ell}_{^{\star, T}}}=T$.  
\end{prro}

{\bf Proof.} \rm Assume that $T^{{\ell}_{^{\star, T}}}=T$, that is
${\ell}_{^{\star, T}}$ is a (semi)star operation of finite type on $T$
(Proposition \ref{pr:3.10} (4)).  In this situation, we have ${\ell}_{^{\star, T}}\le t_T$
and thus $T$ is $({\ell}_{^{\star, T}}, t_T)$--linked to $T$.  By
Proposition \ref{pr:3.10} (3), $T$ is $(\star, {\ell}_{^{\star,
T}})$--linked to $D$.  By transitivity (Lemma \ref{lm:3.1} (b)), we
conclude that $T$ is $t$--linked to $(D, \star)$.  \hfill $\Box$

%%%%%%%%%%%%%%%%%%%%%%%%%%%%%%%%%%
%%%%%%%%%%%%%%%%%%%%%              SECTION 4
%%%%%%%%%%%%%%%%%%%%%%%%%%%%%%%%%%

 \begin{center}
\section{Semistar flatness}
 \end{center}

\hskip .5cm Let $D$ be an integral domain and $T$ be an overring of $D$ and let
$\star$ (respectively, $\star'$) be a semistar operation on $D$ (respectively,
on $T$).  We say that $T$ is {\it $(\star, \star')$--flat over $D$} if,   for  each quasi--$\star'_{f}$--prime
ideal $Q$ of $T$, $(Q\cap D)^{\star_{f}} \neq D^{\star}$ (i.e.  $T$
is $(\star, \star')$--linked to $D$) and, moreover, 
$D_{Q\cap D}= T_{Q}$.  

\vskip -3pt We say that $T$ is {\it $t$--flat over $D$}, if $T$ is
$(t_{D}, t_{T})$--flat over $D$.  Note that, from \cite[Remark 2.3]{KP},
this definition of $t$--flatness coincides with that introduced in
\cite{KP}.  More generally, we say that $T$ is \it $t$-flat over $(D,
\star)$ \rm if $T$ is $(\star, t_{T})$--flat over $D$.  

\smallskip

%REMARK 4.1
\begin{reem} \label{rk:4.1} \rm \bf (a) \rm If $\star := d_{D}$ (respectively,
$\star' := d_{T}$) the identity (semi)star operation on $D$ (respectively,
$T$), then $T$ is $(d_{D}, d_{T})$--flat over $D$ if and only if $T$ is
flat over $D$.

\bf (b) \rm Note that $T$ is $t$-flat over $(D, \star)$ implies $T$ is
$t$--flat over $D$ (for a converse see the following Lemma \ref{lm:4.2}
(e)).  As a matter of fact, for each $Q\in {\cal{M}}(t_{T})$,  $ D_{{Q\cap
D}}=T_{Q}$ and thus, by \cite{KP}, $T$ is a $t$--flat overring of $D$.

\rm \bf (c) \rm Recall that an example given by Fossum \cite[page
32]{Fossum} shows that, even for a Krull domain (hence, in particular,
for a P$v$MD), $t$--flatness does not imply flatness (cf.  also
\cite[Remark 2.12]{KP}).
\end{reem}

\smallskip

\hskip .5cm The proof of the following lemma, in which we collect some
preliminary properties of semistar flatness, is straightforward.

     %LEMMA 4.2
\begin{leem} \label{lm:4.2} Let $T, S$ be two overrings of an integral
domain $D$, with $D \subseteq T \subseteq S$.  \balf \rm\bf \item \sl
Let $D = T$ and $\star', \star''$ be two semistar operations on $T$. 
Then $T$ is $(\star', \star'')$--flat over $T$ if and only if $T$ is
$(\star', \star'')$--linked to $T$.  This happens when $\star'_{f} \leq
\star''_{f}$.

\rm\bf \item \sl Let  $\star$ (respectively, $\star'$,
$\star''$) be a semistar operation on $D$ (respectively, $T$, $S$). 
Assume that $S$ is $(\star', \star'')$--flat over  $T$ and that $T$ is
$(\star, \star')$--flat over  $D$, then $S$ is $(\star, \star'')$--flat over  $D$.

\rm\bf \item \sl Let  $\star$ (respectively, $\star'$,
$\star''$) be a semistar operation on $D$ (respectively, two
semistars operations on $T$).  Assume that $\star'_{f} \leq
\star''_{f}$.  If\ $T$ is $(\star, \star')$--flat over $D$, then $T$ is also
$(\star, \star'')$--flat over $D$.

\rm\bf \item \sl Let $\star$ be a semistar operation on $D$ and let
$\star'$ be a (semi)star operation on $T$ (hence, $\star'_{f}\leq t_{T}$).  If
$T$ is $(\star, \star')$--flat over $D$ then $T$ is $t$--flat over $(D,
\star)$.

%\rm\bf \item \sl Let $\star$ be a semistar operation on $D$.  Then
%$D^\star$ is $(\star, \dot{\star})$--flat over $D$.

\rm \bf \item \sl  Let 
$\star_1$ and $\star_2$ be two semistar operations on $D$ and let $\star'$
be a semistar operation on $T$.
Assume
that $(\star_1)_{f} \leq (\star_2)_{f}$.  If $T$ is $(\star_2, \star')$--flat
over $D$, then $T$ is $(\star_1, \star')$--flat over $D$.  In particular
(cf.  also Remark \ref{rk:4.1} (b)), if $\star$ is a (semi)star operation on $D$ (hence
$\star_{f }\le t_D$), then $T$ is $t$--flat over $(D, \star)$ if and
only if $T$ is $t$--flat over $D$.

\rm \bf \item \sl Let  $\star$ (respectively, $\star'$) be a semistar operation on $D$ 
(respectively, $T$). The overring $T$ is $(\star, \star')$--flat over $D$ if
and only if, for each quasi--$\star'_{f}$--maximal ideal $N$ of $T$, $(N\cap
D)^{\star_{f}} \neq D^{\star}$ and $D_{N\cap D}= T_{N}$.

\rm \bf \item \sl   Let  $\star$ (respectively, $\star'$,
$\star''$) be a semistar operation on $D$ (respectively, $T$, $S$). 
Assume that $S$ is $(\star, \star'')$--flat over $D$ and that each
quasi--$\star'_{f}$--maximal ideal of $T$ is the contraction of a
quasi--$\star''_{f}$--(maximal)ideal of $S$, then $T$ is $(\star,
\star')$--flat over $D$.

% \vskip -3pt In particular, if we take $S:=T^{\star'}$ and
% $\star'':={\dot{\star'}}$ (note that ${\dot{\star'}}$ is a (semi)star
% operation on $T^{\star'}$), then $T$ is $(\star, \star')$--flat over  $D$
% if and only if $T^{\star'}$ is $(\star, {\dot{\star'}} )$--flat over
% $D$. 

\rm \bf \item \sl Let $\star$ (respectively, $\star'$, $\star''$) be a
semistar operation on $D$ (respectively, $T$, $S$).  Assume that $S$ is
$(\star, \star'')$--flat over $D$.  Then $S$ is
$(\star', \star'')$--flat over $T$ if and only if $S$ is
$(\star', \star'')$--linked with $T$.  \hfill $\Box$

\ealf
\end{leem}

\smallskip

%REMARK 4.3
\begin{reem} \label{rk:4.3} \rm \bf (a) \rm When $\star$ is a proper semistar
operation on $D$ (that is $D^{\star}\neq D$), the equivalence of the
second part of statement (e) in the previous lemma fails to be true in
general.  Indeed, if $\star=e_{D}$ then each $t$--flat overring $T$ of
$D$ is not $t$--flat over $(D, e_{D})$, since $T$ is not
$(e_{D},t_{T})$--linked with $D$.  An example in case $\star \neq e_{D}$
is given next.

\vskip -3pt Let $D$ be a Pr\"ufer domain with two prime ideals 
$P\not\subseteq Q$.  Let $T:=D_P$ and consider $\star:=\star_ {\{D_Q\}}$ as a
semistar operation of finite type on $D$.  Then $T$ is $t$--flat over $D$ (since
$T$ is flat over $D$), but $T$ is not $t$--flat over $(D, \star)$. 
Indeed, we have that $M:=PD_P$ is a $t$--ideal of $T$ and $(M\cap
D)^{\star}=P^{\star}=PD_Q=D_Q=D^{\star}$.

\rm \bf (b) \rm Note that, for each semistar operation $\star$ on $D$,
$D^\star$ is $(\star, \dot{\star})$--linked to $D$ (Lemma \ref{lm:3.1}
(e)), but in general $D^\star$ is not $(\star, \dot{\star})$--flat over
$D$.  For instance, if $T$ is a proper non-flat overring of $D$ and if
$\star := \star_{\{T\}}$, then $D^\star = T$, $\dot{\star} =d_{T}$ and
$T$ is not $(\star_{\{T\}}, d_{T})$--flat over $D$.

\bf (c) \rm Let $\{ T_{\lambda}\,\mid\,\, \lambda \in \Lambda \}$ be a
family of overrings of $D$ and let $\ast_{\lambda}$ be a semistar
operation defined on $T_{\lambda}$, for $\lambda \in \Lambda $.  Set
$T:= \cap \{ T_{\lambda}\,\mid\,\, \lambda \in \Lambda \}$ and denote by
$\ast_{\Lambda}$ the semistar operation on $T$ associated to the family
$\{(T_{\lambda}, \ast_{\lambda}) \,\mid \, \lambda \in \Lambda \}$
(Example \ref{ex:2.1} (4)).  If $T_{\lambda}$ is $(\star,
\ast_{\lambda})$--flat over $D$, for each $\lambda \in \Lambda$,  is
$T$ $(\star, \ast_{\Lambda})$--flat over $D$ ?

\vskip -3pt The answer is negative, in general.  For instance, let
$V:={\mathbb C}+M$ be a valuation domain with unbranched  maximal ideal $M$
and let $D:={\mathbb R}+M \subseteq V$.  By \cite[Exercise 5 (a), p. 
340]{GILMER:1972}, the domain $D$ has the $QQR$--property, but it is not
a Pr\"ufer domain.  By \cite[Proposition 2.8]{KP}, there exists an
overring $T $ of $D$ which is not $t$--flat (note that, necessarily, $T
= \cap\{D_P\mid\, P\in \Lambda\}$ for some subset $\Lambda$ of the prime
spectrum of $D$).  Let $\star:=t_D$ and let $\ast_P:=t_{D_P}$, for each
$P\in \Lambda$.  Then, obviously, $D_P$ is ($(\star, \ast_P)$--)flat
over $D$, for each $P\in \Lambda$, but $T$ is not $(\star,
\ast_{\Lambda})$--flat over $D$.  Indeed, we have
$(\ast_{\Lambda})_f\leq t_T$, so if $T$ was $(\star,
\ast_{\Lambda})$--flat over $D$, then $T$ would be $t$--flat over $D$
(Lemma \ref{lm:4.2} (d)).

\end{reem}

\medskip

\hskip .5cm Let $\star$ be a semistar operation on an integral domain $D$
with field of quotients $K$, if $\Sigma$ is a multiplicative system of
ideals of $D$, then we set $\Sigma^\star:=\{I^\star\mid \, I\in \Sigma\}$. 
It is easy to verify that $\Sigma^\star$ is \it a
$\dot{\star}$-multiplicative system of $\dot{\star}$-ideals of $D^\star$
\rm (i.e., if $I^\star, J^\star \in \Sigma^\star$ then $(I^\star\cdot
J^\star)^{\dot{\star }}= (I\cdot J)^\star \in \Sigma^\star$). 

\vskip -3pt If $\Sigma$ is a multiplicative system of ideals of $D$, then:
$$
D^\star_{\Sigma^\star}:= \{z \in K\mid \, zI^\star \subseteq D^\star\,,
\mbox{ for some } I \in \Sigma \}
$$
is an overring of $D^\star$ (and of $
D_{\Sigma}:= \{z \in K \mid \, zI\subseteq D\,, \mbox{ for some } I \in
\Sigma \} $), called \it the generalized                                              ring of fractions of $D^\star$ with
respect to the $\dot{\star}$-multiplicative system $\Sigma^\star$. \rm

\medskip

%PROPOSITION 4.4
\begin{prro} \label{pr:4.4} \sl Let $D$ be an integral domain and $T$ be 
an overring of $D$. Let $\star$ (respectively, $\star'$) be a 
semistar operation on $D$ (respectively, on $T$). The following statements are 
equivalent:
\brom
\rm \bf \item \sl $T$ is $(\star, \star')$--flat over $D$;

\rm \bf \item \sl $T$ is $(\star, \star')$--linked with $D$ and, for each
prime ideal $P$ of $D$, either $(PT)^{\star'_f} = T^{\star'} $ or
$T\subseteq D_{P}$;

\rm \bf \item \sl $T$ is $(\star, \star')$--linked with $D$ and, for each $x\in
T$, $x\neq 0$, $((D:_{D}xD)T)^{\star'_f}=T^{\star'}$;

\rm \bf \item \sl $T$ is $(\star, \star')$--linked with $D$ 
and $T^{\widetilde{\star'}} =\cap \{D_{Q\cap D}\; | \;\; Q \in
{\cal{M}}(\star'_{f})\}$;

\rm \bf \item \sl $T$ is $(\star, \star')$--linked with $D$ and, there exists 
a multiplicative system of ideals $\Sigma$ in $D$ such that
$T^{\widetilde{\star'}}= D^{\tilde{\star}}_{\Sigma^{\tilde{\star}}}$ and
$(IT)^{{\star'_f}}=T^{\star'}$ for each $I\in\Sigma$.

\smallskip \hskip -0.8cm
Moreover, each of the previous statement is a consequence of the
following:

\rm \bf \item \sl  $T$ is $(\star, \star')$--linked with $D$ and, for each
quasi--${\star_{f}}$--prime ideal $P$ of $D$, $T_{D\setminus P}$ is flat
over $D_{P}$.

\erom
\end{prro}

{\bf  Proof.} (i) $ \Rightarrow$ (ii).  Let $P$ be a prime ideal of
$D$.  Assume that $(PT)^{{\star'_f}} \neq T^{{\star'}} $ then there
exists $Q \in {\cal{M}}(\star'_{f})$ such that $PT \subseteq Q$, and so $P
\subseteq Q\cap D$.  Therefore, by the assumption,  $D_{P}\supseteq D_{Q\cap D}=
T_{Q}\supseteq T$.

(ii) $\Rightarrow$ (iii).  Let $0\neq x\in T$. 
Assume that $((D:_{D}xD)T)^{{\star'_f}} \neq T^{{\star'}}$, then there exists
$Q \in {\cal{M}}(\star'_{f})$ such that $(D:_{D}xD)T \subseteq Q$.  We have
$(D:_{D}xD) \subseteq Q\cap D=:P$ and $(PT)^{{\star'_f}} \neq T^{{\star'}} $.
Hence, by assumption, $T\subseteq D_{P}$.  Write $x=\frac{d}{s}$, for
some $d\in D$ and $s\in D\setminus P$.  Then $s\in (D:_{D}xD)\subseteq P$,
which is impossible.

(iii) $\Rightarrow$ (iv).  By the definition of $\widetilde{\star'}$ we have that
$T^{\widetilde{\star'}}=\cap\{T_Q\mid\, Q \in {\cal{M}}(\star'_{f})\}$,
and hence $\cap \{D_{Q\cap D}\mid \, Q \in
{\cal{M}}(\star'_{f})\}\subseteq T^{\widetilde{\star'}}$.  For the
reverse inclusion, let $x\in T$, $x\neq 0$,  then
$((D:_{D}xD)T)^{{\star'_f}}=T^{\star'}$.  Let $Q \in
{\cal{M}}(\star'_{f})$.  Then $(D:_{D}xD)T \not\subseteq Q$, that is
$(D:_{D}xD)\not\subseteq Q\cap D$.  So $x\in D_{Q\cap D}$.  Thus
$T\subseteq D_{Q\cap D}$, and hence $T_Q =D_{Q\cap D}$.  Therefore
$T^{\widetilde{\star'}}\subseteq D_{Q\cap D}$ for each $Q \in
{\cal{M}}(\star'_{f})$ and so we conclude that $T^{\widetilde{\star'}}
=\cap \{D_{Q\cap D}\mid \, Q \in {\cal{M}}(\star'_{f})\}$.

(iv) $\Rightarrow$ (i).  Let $Q \in
{\cal{M}}(\star'_{f})$.  Then $T\subseteq T^{\widetilde{\star'}}
\subseteq D_{Q\cap D}$.  Hence $T_Q\subseteq D_{Q\cap D}$.  The reverse
inclusion is trivial.

(ii) $\Rightarrow$ (v).  Let
$\Sigma :=\{I  \hbox{ nonzero ideal of }\ D \mid \,
(IT)^{{\star'_f}}=T^{\star'}\}$.  The set $\Sigma$ is a multiplicative
system of ideals of $D$.  Hence
$\Sigma^{\tilde{\star}}=\{I^{\tilde{\star}}\mid \, I\in \Sigma\}$ is a
${\dot{\tilde{\star}}}$--multiplicative system of
${\dot{\tilde{\star}}}$--ideals of $D^{\tilde{\star}}$.  Let $x\in
D^{\tilde{\star}}_{\Sigma^{\tilde{\star}}}$.  Then $xI\subseteq
xI^{\tilde{\star}} \subseteq D^{\tilde{\star}}$, for some $I\in \Sigma$. 
Since $D^{\tilde{\star}}\subseteq T^{\widetilde{\star'}}$ (Remark
\ref{rk:3.3} (b)), then $xIT\subseteq T^{\widetilde{\star'}} $, and
hence $x(IT)^{\widetilde{\star'}}\subseteq T^{\widetilde{\star'}} $.  On
the other hand, since $(IT)^{{\star'_f}}=T^{\star'}$, then necessarily
$(IT)^{\widetilde{\star'}}= T^{\widetilde{\star'}} $.  Hence
$xT^{\widetilde{\star'}}\subseteq T^{\widetilde{\star'}} $ and so $x\in
T^{\widetilde{\star'}} $.  Therefore
$D^{\tilde{\star}}_{\Sigma^{\tilde{\star}}}\subseteq
T^{\widetilde{\star'}}$.

\vskip -4pt For the opposite inclusion, let $0\neq x\in
T^{\widetilde{\star'}} $.  Set $I:=(D:_DxD)$.  We claim that
$(IT)^{{\star'_f}}=T^{\star'}$ (i.e. $I \in \Sigma$).  Otherwise, as in
the proof of (ii) $\Rightarrow$ (iii), there exists $Q\in
{\cal{M}}(\star'_{f})$ such that $I\subseteq Q\cap D$ and $T \subseteq
D_{Q\cap D}$.  Hence $T^{\widetilde{\star'}}\subseteq T_Q\subseteq
D_{Q\cap D}$.  Write $x=\frac{d}{s}$ for some $d\in D$ and $s\in
D\setminus (Q\cap D)$.  Therefore $s\in (D:_DxD)\subseteq Q\cap D$,
which is impossible.

\vskip -4pt Finally, in general, we have
$xI^{\widetilde{\star}}=(xI)^{\widetilde{\star}} \subseteq 
(x(D:_{K}xD))^{\widetilde{\star}}= D^{\widetilde{\star}}$.  So $x\in
D^{\tilde{\star}}_{\Sigma^{\tilde{\star}}}$ (i.e.
$T^{\widetilde{\star'}}\subseteq D^{\tilde{\star}}_{\Sigma^{\tilde{\star}}}$),
hence we conclude that $T^{\widetilde{\star'}}=
D^{\tilde{\star}}_{\Sigma^{\tilde{\star}}}$.

(v) $\Rightarrow$ (iv). The inclusion  $\cap \{D_{Q\cap D}\; | \;\; Q 
\in {\cal{M}}(\star'_{f})\}\subseteq T^{\widetilde{\star'}}$ is 
clear. Now, let $x\in T^{\widetilde{\star'}}= D^{\tilde{\star}}_{\Sigma^{\tilde{\star}}}$. 
Then there exists a nonzero ideal $I\in \Sigma$ such that $xI\subseteq
D^{\widetilde{\star}}$.  Let $Q\in {\cal{M}}(\star'_{f})$.  Since $I \in
\Sigma$ then, by assumption, $(IT)^{{\star'_f}}=T^{\star'}$ and, thus,
$I\not\subseteq Q\cap D$.  Let $s\in I\setminus(Q\cap D)$, then $sx\in
D^{\widetilde{\star}}$.  On the other hand, since $(Q\cap
D)^{\star_f}\neq D^\star$ (Proposition \ref{pr:3.2}), there exists $M\in
{\cal{M}}(\star_{f})$ such that $Q\cap D\subseteq M$.  Therefore we have
that $D^{\widetilde{\star}}\subseteq D_M\subseteq D_{Q\cap D}$ and so
$sx\in D_{Q\cap D}$, thus $x\in D_{Q\cap D}$.  Hence we conclude that
$T^{\widetilde{\star'}} \subseteq\cap \{D_{Q\cap D}\; | \;\; Q \in
{\cal{M}}(\star'_{f})\}$.

(vi) $\Rightarrow$(i).  Let $Q$
be a quasi--$\star'_{f}$--prime ideal, and let $P$ be a
quasi--$\star_f$--prime of $D$ such that $(Q\cap D)^{\star_f}\subseteq
P$ (Proposition \ref{pr:3.2}).  Since  $QT_{D\setminus
P}$ is a prime ideal of $T_{D\setminus P}$ such that $QT_{D\setminus P}
\cap D_{P} =(Q\cap D)D_{P}$ and, by assumption, $T_{D\setminus P}$ is
flat over $D_{P}$, then we conclude that $D_{Q\cap D}=(D_{P})_{(Q\cap
D)D_{P}}= (T_{D\setminus P})_{QT_{D\setminus P} }=T_{Q}$. 

\vskip -4pt \hfill $\Box$

\smallskip

%THEOREM 4.5
\begin{thee} \label{th:4.5} \sl Let $D$ be an integral domain and $T$ be 
an overring of $D$. Let $\star$ (respectively, $\star'$) be a 
semistar operation on $D$ (respectively, on $T$). The following statements are 
equivalent:
\brom
\bf \item \sl  $T$ is $(\star, \star')$--flat over $D$;

\bf \item \sl  
$\Na(T, \star')$ is a flat overring of $\Na(D, \star)$;

\bf \item \sl $T$ is $({\widetilde{\star}}, {\widetilde{\star'}})$--flat
over $D$;

\bf \item \sl  
$T^{\widetilde{\star'}}$ is a $(\dot{\tilde{\star}},
\dot{\widetilde{\star'}})$--flat overring of $D^{\tilde{\star}}$.  \erom
\end{thee}

\bf Proof.  \rm  
Since $\Na(D, \star)=\Na(D, \widetilde{\star}) =\Na(D^{\widetilde{\star}},
\dot{\widetilde{\star}})$ and, similarly, $\Na(T, \star')=\Na(T,
\widetilde{\star'}) $ $=\Na(T^{\widetilde{\star'}}, \dot{\widetilde{\star'}})$
(Lemma \ref{lm:2.4} (i)), it suffices to show that (i) $\Leftrightarrow$ (ii).

(i) $\Rightarrow$ (ii). Since $T$ is  $(\star, \star')$--linked to 
$D$, then $\Na(D, \star)\subseteq \Na(T, \star')$, by Theorem \ref{th:3.7}.  Now,
let $N$ be a maximal ideal of $\Na(T, \star')$.  Then $N =Q\Na(T,
\star')=QT_Q(X)\cap \Na(T, \star')$, for some
$Q\in{\cal{M}}(\star'_{f})$ (cf.  also Lemma \ref{lm:2.4} (e)), and
$\Na(T, \star')_N =\Na(T, \star')_{Q\Na(T, \star')}=T_Q(X)=D_{Q\cap
D}(X)$, because of Corollary \ref{cor:2.5} and, by assumption, $D_{Q\cap
D}=T_Q$ .  On the other hand, by semi\-star linkedness, $(Q\cap
D)^{\star_f}\neq D^\star$ (Proposition \ref{pr:3.2}) then we have that
$\Na(D, \star)_{(Q\cap D)\Na(D, \star)}$ $= D_{Q\cap D}(X)$ (Corollary
\ref{cor:2.5}).  One can easily check that $N \cap \Na(D, \star)= (Q\cap
D)\Na(D, \star)$.  Therefore $\Na(T, \star')_N= \Na(D, \star)_{N \cap
\Na(D, \star)}$, as desired.

(ii) $\Rightarrow$ (i).  Since $\Na(D, \star)\subseteq \Na(T, \star')$,
then $T$ is $(\star, \star')$--linked to $D$ (Theo\-rem \ref{th:3.7}). 
Let $Q$ be a quasi--${\star'_{f}}$--maximal ideal of $T$ and set
$N:=Q\Na(T, \star')$, then $\Na(T,
\star')_{N}= \Na(T, \star')_{Q\Na(T, \star')}= T_{Q}(X)$ (Corollary \ref{cor:2.5}).
On the other hand, by flatness, we have $\Na(T, \star')_{N}=
\Na(D, \star)_{N \cap \Na(D, \star)}=\Na(D, \star)_{(Q\cap D)\Na(D,
\star)}$.  Since, by semistar linkedness, $(Q\cap D)^{\star_f}\neq
D^\star$ (Proposition \ref{pr:3.2}), then we have that $\Na(D,
\star)_{(Q\cap D)\Na(D, \star)}=D_{Q\cap D}(X)$ (Corollary \ref{cor:2.5}).  Therefore
$T_Q(X)=D_{Q\cap D}(X)$ and so $T_Q=D_{Q\cap D}$.  Hence $T$ is $(\star,
\star')$--flat over $D$.  \hfill $\Box$ \bigskip

\hskip .5cm  The following result sheds new light on the statement (vi) of Proposition
\ref{pr:4.4}.

% PROPOSITION 4.6
\begin{prro} \label{pr:4.8} Let $D$ be an integral domain and $T$ be an
overring of $D$.  Let $\star$  be a semistar
operation on $D$ and let $\ell :=\ell_{\star, T}$
be the semistar operation on $T$ introduced in Section 3. The following
statements are equivalent: 

\brom 
% i
 \bf \item \sl $T$ is  $(\star, \ell)$--flat over $D$;
 
% ii
\bf \item \sl for each prime ideal $P$
of $D$, either $(PT)^{\ell} = T^{\ell} $ or $T\subseteq D_{P}$;

% iii
\bf \item \sl 
 for each $x\in T$,\ $x \neq 0$,\
 $((D:_{D} xD)T)^{\ell}= T^{\ell}$;
 
% iv
\bf \item \sl $T^{\ell} = \cap \{ D_{{N \cap D}} \mid \, N \, \mbox{\sl  is a
prime ideal of } \, T\,, \; \mbox{\sl  maximal with the property   }$ \newline
${ }$ \hskip 2.2cm $(N \cap D)^{\star_{f}} \neq D^\star \}$;

% v
\bf \item \sl for each prime ideal $Q$
of $T$ such that $(Q\cap D)^{\star_{f}} \neq D^\star$, then  $D_ {Q\cap D}= T_Q$;

%vi
\bf \item \sl 
 for each prime ideal $N$ of $T$, \,
maximal with respect to the property \, $(N\cap D)^{\star_f}\neq D^\star$,
then $D_ {N\cap D}= T_N$;

%vii
\bf \item \sl  for each
quasi--${\star_{f}}$--prime ideal $P$ of $D$, $T_{D\setminus P}$ is flat
over $D_{P}$;

% viii
\bf \item \sl 
for each nonzero finitely
generated fractional ideal $F$ of $D$,
$((D:_{K}F)T)^{\ell}=(T^\ell:_{K}FT)$.

\erom
\end{prro}

{\bf Proof.} Note that the set of quasi--$\ell$--prime (respectively,
quasi--$\ell$--maximal) ideals of $T$ coincides with the set of prime
ideals $Q$ of $T$ such that $(Q\cap D)^{\star_{f}} \neq D^\star$
(respectively, the set of prime ideals $N$ of $T$, maximal with the
property $(N \cap D)^{\star_{f}} \neq D^\star)$.  Therefore, the
statements (i) -- (vi) are equivalent by Proposition \ref{pr:4.4} and
Proposition \ref{pr:3.10} (3).

(v) $\Rightarrow$ (vii).  Let $P$ be a quasi--$\star_{f}$--prime ideal
of $D$.  Let $N$ be a maximal ideal of $T_{D\setminus P}$.  Then ${N\cap
D}\subseteq P$ and, hence, $((N\cap T)\cap D)^{\star_f}\neq D^\star$. 
So $D_{N\cap D} = T_{N\cap T}$.  On the other hand, we have
$(T_{D\setminus P})_N=T_{N\cap T}$, and $(D_P)_{N\cap D_P}=(D_P)_{(N\cap
D)D_P}=D_ {N\cap D}$.  Hence $(D_P)_{N\cap D_P}=(T_{D\setminus P})_N$,
as desired.

(vii) $\Rightarrow$ (viii).  We have
$((D:_{K}F)T)^{\ell}=\cap\{(D:_{K}F)T_{D\setminus P}\mid \, P\ \mbox{ is
a quasi--$ \star_f$--}$ prime of $\ D\}$.  As $T_{D\setminus P}$ is
$D_P$--flat (hence, $T_{D\setminus P}$ is also $D$--flat) and $F$ is
finitely generated, then $(D:_{K}F)T_{D\setminus P}= (T_{D\setminus
P}:_{K}FT_{D\setminus P})=(T:_{K}FT)T_{D\setminus P}$, for each quasi--$
\star_f$--\-prime $P$ of $D$.  Hence
$((D:_{K}F)T)^{\ell}=\cap\{(T:_{K}FT)T_{D\setminus P}\, \mid\; P\ \mbox{
is a quasi--$ \star_f$--prime of }\ D\}$ =
$(T:_{K}FT)^{\ell}=(T^{\ell}:_{K}FT)$ (since $\ell$ is stable; Example
\ref{ex:2.1} (5) and Proposition \ref{pr:3.10} (1)).

(viii) $\Rightarrow$ (iii).  Take $F:= D + xD$.    \hfill $\Box$

\bigskip
\hskip .5cm It is well-known that a domain with all its overrings
flat (or, equivalently, with all its overrings $t$--flat) coincides with
a Pr\"ufer domain (cf.  \cite [Theorem 4]{R}, \cite[Proposition
2.8]{KP}).  The following proposition deals with a similar question in
the semistar case.

\smallskip

% THEOREM 4.7
\begin{thee} \label{th:4.9}
\sl  Let $D$ be an integral domain and $\star$ a semistar operation on
$D$.  The following statements are equivalent:
\brom
\bf \item \sl For each overring $T$ of $D$ and for each semistar operation
$\star'$ on $T$, $T$ is $(\star, \star')$--flat over $D$;

\bf \item \sl  Each
overring $T$ of $D$ is $(\star, d_T)$--flat over $D$;

\bf \item \sl Each overring $T$ of $D$ is $t$--flat over
$(D,\star)$;

\bf \item \sl $D$ is a Pr\"ufer domain in which each maximal
ideal is a quasi--$\star_f$--maximal ideal.
\erom
\end{thee}

{\bf Proof.}  (i) $\Rightarrow$ (ii) is obvious. 

(ii) $\Rightarrow$ (iii) is a consequence of $d_T\le t_T$ (Lemma
\ref{lm:4.2} (c)).

(iii) $\Rightarrow$ (iv).  Since semistar flatness implies
semistar linkedness, then, by Theorem \ref{th:3.8}, each maximal ideal is a
quasi--$\star_f$--maximal ideal.  On the other hand, since an overring
$t$--flat over $(D,\star)$ is also $t$--flat over $D$ (Remark
\ref{rk:4.1} (b)), then each overring of $D$ is $t$--flat over $D$. 
Hence, by \cite[Proposition 2.8]{KP}, $D$ is a Pr\"ufer domain.

(iv) $\Rightarrow$ (i).  Let $T$ be an overring of $D$ and $\star'$ a
semistar operation on $T$.  Let $Q$ be a quasi--$\star'_f$--prime ideal
of $T$.  Then $Q\cap D$ is contained in a maximal ideal of $D$ which is,
by assumption, a quasi--$\star_f$--maximal ideal of $D$.  Therefore $(Q\cap
D)^{\star_f}\neq D^\star$, and so $T$ is $(\star, \star')$--linked with $D$. The
equality $T_Q=D_{Q\cap D}$ is a consequence of the fact that $T$ is an
overring of the Pr\"ufer domain $D$ \cite[Theorem 26.1]{GILMER:1972}.  \hfill
$\Box$

\bigskip

%%%%%%%%%%%%%%%%%%%%%%%%%%%%%%%%%%
%%%%%%%%%%%%%%%%%%%%%              SECTION 5
%%%%%%%%%%%%%%%%%%%%%%%%%%%%%%%%%%

 %\begin{center}
\section{Pr\"ufer semistar multiplication domains}
%\end{center}

\hskip .5cm As an application of the previous sections, our goal is to
give  new characte\-ri\-zations of Pr\"ufer semistar multiplication  domains,
in terms of semistar linked overrings and semistar flatness.

\hskip .5cm Let $D$ be an integral domain and $\star$ a semistar
operation on $D$.  Recall that $D$ is \it a P$\star$MD
(Pr\"ufer
$\star$--multiplication domain)\rm , if each $F\in
\boldsymbol{f}(D)$ is $\star_f$-invertible (i.e.,
$(FF^{-1})^{\star_f}=D^\star$).

\vskip -3pt The notion of P$\star$MD is a generalization of the notion of
Pr\"ufer
$v$--multiplication domain (cf.  \cite[page 427]{GILMER:1972},
\cite{Griffin:1967}, \cite{Mott/Zafrullah:1981}) and so, in particular, of Pr\"ufer
domain.  When $\star=d$ (where $d$ is the identity (semi)star operation
on $D$) the P$d$MDs are just the Pr\"ufer domains.  If $\star=v$ (where
$v$ is the $v$-(semi)star operation on $D$), we obtain the notion of
P$v$MD.

% REMARK 5.1

\begin{reem}\label {rk:5.1} \rm \bf (a) \rm The notions of
P$\star$MD and P$\star_f$MD coicide.  In particular, a P$v$MD coincides
with a P$t$MD.

\bf (b) \rm Let $\star_1$ and $\star_2$ be two semistar operations on $D$
such that $\star_1\le\star_2$.  If $D$ is a P$\star_1$MD, then $D$ is also
a P$\star_2$MD. In particular, if $\star$ is a (semi)star operation on
$D$, and hence $\star \le v$ \cite[Theorem 34.1 (4)]{GILMER:1972}, then a P$\star$MD is
a P$v$MD. Also, since $d \le\star$ for any semistar operation $\star$,
then a Pr\"ufer domain is a P$\star$MD for any arbitrary semistar
operation $\star$ on $D$.

\bf (c) \rm In the semistar case (i.e. if $\star$ is a proper semistar
operation), a P$\star$MD is not necessarily integrally closed \cite[Example 3.10]{FJS}.
\end{reem}

\hskip .5cm We recall some of the characterizations of P$\star$MDs
proved in \cite{FJS}:

\smallskip

% THEOREM 5.2
\begin{thee}\label {th:5.2} \rm \cite [Theorem 3.1, Remark 3.2]{FJS} \sl Let
$D$ be an integral domain and $\star$ a semistar operation on $D$.  The
following statements are equivalent: \brom \bf \item \sl $D$ is a
P$\star$MD;

\bf \item \sl $D_Q$ is a valuation domain, for each $Q\in {\cal
M}(\star_f)$;

\bf \item \sl $\Na(D, \star)$ is a Pr\"ufer domain;\

\bf \item \sl  $D$ is a
P$\tilde{\star}$MD.
\erom
Moreover, if $D$ is a P$\star$MD, then $\tilde{\star}=\star_f$.  \hfill $\Box$
\end{thee}

\smallskip

\hskip .5cm The following theorem is ``a semistar version'' of a  characterization of the Pr\"ufer
domains proved by E. Davis  \cite[Theorem 1]{Davis}.  It generalizes
properly \cite[Theorem 2.10]{DHLZ1}, stated in the case of
$t$-operations (cf.  also \cite [Theorem 5.1]{Mott/Zafrullah:1981} and
\cite[Corollary 3.9]{Kang:1989}).

\vskip -3pt Recall that an integral domain $D$, with field of quotients
$K$, is \it seminormal \rm if, whenever $x \in K$ satisfies $x^{2}, x^{3} \in
D$, then $x\in D$, \cite{Gilmer/Heitmann}.
 
 \smallskip

% THEOREM 5.3
\begin{thee}\label {th:5.3} \sl Let $D$ be an integral domain , $T$ an
overring of $D$,  $\star$ a semistar operation on $D$ and
let $\ell :=\ell_{\star, T}$ be the semistar operation on $T$ introduced
in Section 3.  The following statements are equivalent: \brom \bf \item
\sl For each overring $T$ and for each semistar operation $\star'$ such that
$T$ is $(\star, \star')$--linked to $D$, $T^{\widetilde{\star'}}$ is
integrally closed.

\bf \item \sl  For each overring
$T$ of $D$, $T^{{\ell}_{^{\star, T}}}$ is integrally closed.

\bf \item \sl 
Each overring $T$, $t$--linked to $(D, \star)$, is integrally closed.

\bf \item \sl  Each overring $T$, $(\star, d_T)$-linked to $D$, is integrally
closed.

\bf \item \sl $D^{\widetilde{\star}}$ is integrally closed and, for
each overring $T$ and for each semistar operation $\star'$ on $T$ such that
$T$ is $(\star, \star')$--linked to $D$, $T^{\widetilde{\star'}}$ is
seminormal.

\bf \item \sl  $D^{\widetilde{\star}}$ is integrally closed and
each overring $T$, $t$--linked to $(D, \star)$, is seminormal.

\bf \item \sl 
$D^{\widetilde{\star}}$ is integrally closed and each overring $T$,
$(\star, d_T)$-linked to $D$, is seminormal.

\bf \item \sl  $D$ is a
P$\star$MD. 
\erom
\end{thee}

\bf Proof.  \rm (i) $\Rightarrow$ (ii).  It follows from Proposition
\ref{pr:3.10} (3) and (4), by taking $\star'=\ell_{\star, T}$.

(ii) $\Rightarrow$ (iii) follows from Proposition \ref{pr:3.13}.

(iii) $\Rightarrow$ (iv).  Obvious since $d_T\le t_T$ (Lemma
\ref{lm:3.1} (c)).

(iv) $\Rightarrow$ (v).  Let $(T, \star')$ be such that $T$ is $(\star,
\star')$--linked to $D$.  Let $P$ be a quasi--$\star_f$--prime ideal of
$D$.  By Example \ref{ex:3.4} (1), $T_{D\setminus P}$ is $(\star,
d_{T,P})$--linked to $D$.  Hence, by assumption,  $T_{D\setminus P}$ is
integrally closed.  In particular (for $(T, \star') = (D, \star)$), $D_P$ is
integrally closed, and hence $D^{\widetilde{\star}}$ is integrally
closed.  On the other hand, if $Q$ is a quasi--$\star'_f$--prime ideal
of $T$, there exists $P$ a quasi--$\star_f$--prime ideal of $D$ such
that $Q\cap D\subseteq P$ (Proposition \ref{pr:3.2}).  Hence
$T_{D\setminus P}\subseteq T_Q$ and so $T_Q$ is integrally closed, since
$T_{D\setminus P}$ is.  Therefore, $T^{\widetilde{\star'}}$ is
integrally closed; in particular, $T^{\widetilde{\star'}}$ is
seminormal.

(v) $\Rightarrow$ (vi) is obvious and (vi) $\Rightarrow$ (vii) is a
consequence of $d_T\le t_T$ (Lemma \ref{lm:3.1} (c)).

(vii) $\Rightarrow$ (viii).  We want to
show that, for each quasi--$\star_{f}$--maximal ideal $P$ of $D$, $D_{P}$ is a
valuation domain (Theorem \ref{th:5.2}), i.e., if $x $ is a nonzero
element of the quotient field $K$ of $D$, then either $x$ or $x^{-1}$ is in
$D_{P}$ .  Note that, from the assumption, it follows that $D_{P}=
D^{\tilde{\star}}_{PD_{P}\cap D^{\tilde{\star}}}$ is integrally closed. 
If we set $T:=D[x^{2}, x^{3}]$ then (by Example \ref{ex:3.4} (1)) $T_{D
\setminus P} =D[x^{2}, x^{3}]_{D \setminus P}=D_{P}[x^{2}, x^{3}]$ is
$(\star, d_{T, P})$--linked to $D$ thus, by assumption, $D_{P}[x^{2},
x^{3}]$ is seminormal, i.e., $x \in D_{P}[x^{2}, x^{3}]$.  Hence $x$ is
the root of some polynomial $f $ with coefficients in $D_{P}$ and with
the coefficient of the linear term equal to 1.  This implies that either
$x$ or $x^{-1}$ is in $D_{P}$, by \cite[Theorem 67]{Kaplansky:1970}.

(viii) $\Rightarrow$ (i).  Let $T$ be an overring of $D$ $(\star,
\star')$--linked to $D$.  For each quasi--$\star'_{f}$--maximal ideal
$N$ of $T$, let $P$ be a quasi--$\star_{f}$--maximal ideal $P$ of $D$,
such that $N\cap D \subseteq P$ (Proposition \ref{pr:3.2}), thus $D_{P} \subseteq
T_{D\setminus P}\subseteq T_{N}$.  Since $D$ is a P$\star$MD, then
$D_{P}$ is a valuation domain, hence $T_{N}$ is also a valuation domain
and so $T^{\widetilde{\star'}}= \cap \{T_{N}\mid \, N \in {\cal
M}(\star'_{f})\}$ is integrally closed.  \hfill $\Box$

\smallskip
\hskip .5cm The following result generalizes \cite[Theorem
5.1]{Mott/Zafrullah:1981} (cf.  also \cite[Corollary 3.9]{Kang:1989}).

% COROLLARY 5.4
\begin{coor}\label {cor:5.4} \sl  Let $D$ be an integral domain and $T$ be an 
overring of $D$. Let $\star$ 
(respectively, $\star'$) be a semistar operation on $D$ 
(respectively, on 
$T$).
Assume that $D$ is a P$\star$MD and that $T$ is $(\star,
\star')$--linked to $D$, then $T$ is a P$\star'$MD.
\end{coor}

{\bf Proof.} If $S$ is
an overring of $T$ and $\star''$  a semistar operation on $S$ such that
$S$ is $(\star', \star'')$--linked to $T$, then $S$ is $(\star,
\star'')$--linked to $D$ (Lemma \ref{lm:3.1} (b)). By Theorem 
\ref{th:5.3} ((viii) $\Rightarrow$ (i))  $S^{\widetilde{\star''}}$ is
integrally closed.  The conclusion follows from Theorem 
\ref{th:5.3} ((i) $\Rightarrow$ (viii)).  \hfill $\Box$

\smallskip

% COROLLARY 5.5
\begin{coor}\label{cor:5.5} \sl Let $D$ be P$\star$MD for some semistar
operation $\star$ on $D$.  Then: \balf \bf \item \sl For each overring
$T$ of $D$, $T$ is a P$\dot{\star}^T$MD.

\bf \item \sl  Each $t$--linked overring to $(D, \star)$ is
a P$v$MD. In particular, $D^{[\star]}$ is a P$v$MD and if, moreover,
$(D : D^{\star})\neq 0$, then the complete integral closure $\widetilde{D}$ of
$D$ is a P$v$MD. \ealf
\end{coor}

{\bf Proof.} \bf (a) \rm follows from Corollary \ref {cor:5.4} and Lemma
\ref{lm:3.1} (e).  The first statement in \bf (b) \rm is a particular
case of Corollary \ref {cor:5.4}; the remaining part is a consequence of
the first part and of Corollary \ref{cor:3.6}.  \hfill $\Box$
\medskip

\hskip .5cm Note that Corollary \ref{cor:5.5} (b) generalizes the fact
that the pseudo--integral closure, $D^{[v]}$, of a P$v$MD, $D$, is still
a P$v$MD \cite[Proposition 1.3]{Anderson/Houston/Zafrullah:1991}.

\medskip
%REMARK 5.6
\begin{reem} \label{rk:5.6} \rm The integral closure $D'$ of an integral  domain
$D$ is not in ge\-ne\-ral $t$-linked over $D$ \cite[Example 4.1]{DHLRZ1}.  But, each
domain $D$ has a smallest integrally closed $t$-linked overring, namely
${D'}^{\hskip 1pt \ell_{t_{_{\!\mbox{\tiny {\it D}}}}, D'}}
=\cap\{{D'}_{D\setminus P}\mid\, P \ \mbox{ is a } t-\mbox{\-prime ideal
of } D \} $ \cite[Proposition 2.13 (b)]{DHLZ1}.

\vskip -3pt In the semistar case, $D'$ is always $(\star, {\ell}_{^{{\star},
D'}})$--linked to $D$, for any semistar o\-pe\-ra\-tion $\star$ on $D$
(Proposition \ref{pr:3.10} (3)).  Also note that ${\ell}_{^{{\star}, D'}}$ is the unique
minimal semistar operation in the set of semistar operations $\star'_f$,
where $\star'$ is a semistar operation on $D'$ such that $D'$ is
$(\star, \star')$--linked to $D$ (Proposition \ref{pr:3.10} (5)).  Therefore,
$D'$ is $t$--linked over $D$ if and only if ${\ell}_{t_{_{D}}, D'} \le
t_{D'}$ (Lemma \ref{lm:3.1} (c)) or, equivalently, if and only if ${\ell}_{t_{_{\!D}},
D'}$ is a (semi)star operation on $D'$ (i.e. $D' ={D'}^{\hskip 1pt
\ell_{t_{_{\!D}}, D'}}$).
\end{reem}

\medskip

\hskip .5 cm The next theorem of characterization of P$\star$MDs
 is a ``semistar analogue'' of Richman's flat-theoretic theorem of
 characterization of Pr\"ufer domains \cite [Theorem 4]{R}.  A special
 case of the following result, concerning the $t$-operations, was obtained
 in \cite[Proposition 2.10]{KP}.

\smallskip
%THEOREM 5.7 
\begin{thee} \label{th:5.7} \sl Let $D$ be an integral domain, $\star$ a
semistar operation on $D$, $T$ an overring of $D$ and let $\ell :=\ell_{\star, T}$
be the semistar operation on $T$ introduced in Section 3.  The following
statements are equivalent: \brom 
\bf \item \sl $D$ is a P$\star$MD.

\bf \item \sl For each overring $T$ of $D$ and for each semistar operation
$\star'$ such that $T$ is $(\star, \star')$--linked to $D$, $T$ is
$(\star, \star')$--flat 
over $D$.

\bf \item \sl  For each overring $T$ of $D$, $T$ is
$(\star, {\ell}_{^{\star, T}})$-flat over $D$.

\bf \item \sl For each overring
$T$ of $D$, $t$--linked to $(D, \star)$, $T$ is $t$--flat  over $(D, \star)$.

\bf \item \sl For each overring $T$ of $D$ such that $T$ is $(\star, d_T)$-linked to
$D$, $T$ is $(\star, d_T)$--flat over $D$.

\erom
\end{thee}
{\bf Proof.} (i) $\Rightarrow$ (ii).  Let $T$ be an overring and
$\star'$ a semistar operation on $T$ such that $T$ is $(\star,
\star')$--linked to $D$.  Let $Q$ be   a
quasi--$\star'_f$--prime of $T$ such that $(Q\cap D)^{\star_f}\neq D^\star$. 
Then $Q\cap D\subseteq P$ for some quasi--$\star_{f}$--maximal ideal $P$
of $D$.  Thus $D_{P} \subseteq D_{Q\cap D}\subseteq T_{Q}$.  Since $D$
is a P$\star$MD, then $D_{P}$ is a valuation domain (Theorem \ref{th:5.2}), hence
$T_{Q}$ is also a valuation domain and $T_Q=D_{Q\cap D}$.  Hence $T$ is
$(\star,
\star')$--flat over $D$.

(ii) $\Rightarrow$ (iii) is a trivial consequence of Proposition
\ref{pr:3.10} (3).

(iii) $\Rightarrow$ (ii).  Let $T$ be an overring and $\star'$ a
semistar operation on $T$ such that $T$ is $(\star, \star')$--linked to
$D$.  Then ${\ell}_{^{\star, T}}\le\star'_f$ (Proposition \ref{pr:3.10}
(5)).  Hence $T$ is $(\star, \star')$--flat over $D$ (Lemma \ref{lm:4.2} (c)).

(ii) $\Rightarrow$ (iv) is obvious.

(iv) $\Rightarrow$ (v).  Let $T$ be an overring $(\star, d_T)$-linked to
$D$, and let $Q$ be a prime ideal of $T$.  We have $(Q\cap
D)^{\star_f}\neq D^\star$ (Proposition \ref{pr:3.2}).  Let $P$ be a
quasi--$\star_{f}$--maximal ideal of $D$ such that $Q\cap D\subseteq P$,
thus $D_{P} \subseteq T_{D\setminus P}\subseteq T_{Q}$.  Let $(V, M)$
be a valuation overring of $D$ such that $M\cap D=P$.  Then $D\subseteq
V_{D\setminus P}=V$ is $t$--linked with $(D, \star)$ (Example
\ref{ex:3.4} (1)), and hence $V$ is $t$--flat over $(D, \star)$, by
assumption.  So $V=D_{M\cap D}=D_P$.  Therefore $T_Q \ (\supseteq
D_{P})$ is also a valuation domain and $T_Q=D_{Q\cap D}$, thus $T$ is
$(\star, d_T)$-flat over $D$.

(v) $\Rightarrow$ (i).  Let $P$ be a quasi--$\star_{f}$--prime
ideal of $D$.  Let $T$ be an overring of $D_P$ (and hence of $D$).  Note
that, in this situation, $T=T_{D\setminus P}$.  Hence $T$ is ${(\star,
d_{T})}$--linked to $D$ (Example \ref{ex:3.4} (1)).  So $T$ is ${(\star,
d_{T})}$--flat over $D$, by assumption.  Therefore, if $N$ is a maximal
ideal of $T$, then $T_N=D_{N\cap D}$.  Hence $T_N=(D_P)_{(N\cap
D){D_P}}=(D_P)_{N\cap D_P}$ (since $N\cap D\subseteq P$).  That is, $T$
is $D_P$--flat.  By a result proved by Richman \cite[Lemma 4]{R}, we
deduce that $D_P$ is a valuation domain.  Hence $D$ is a P$\star$MD
(Theorem \ref{th:5.2}).  \hfill $\Box$
\bigskip

%\begin{center}
    \bf Acknowledgement.  \rm \
    %\rm \end{center}
    Both authors are grateful to F. Halter-Koch for providing them with a
    copy of his recent pre-print \it Characterization of Pr\"ufer
    multiplication monoids and domains by means of spectral module theory
    \rm, presented  at the Algebra Conference (Venezia, June
    2002).  Using the language of monoids and module systems, Halter-Koch's
    work sheds further light on some of the themes discussed in the present
    paper.

\end{document}